\documentclass[MSNbibl,number,citesort,seceqn,dvips]{arxbj}
\usepackage{upgreek}
\usepackage{graphicx}

\aid{0}
\volume{19}
\issue{5B}
\pubyear{2013}
\firstpage{2455}
\lastpage{2493}
\doi{10.3150/12-BEJ459} 

\makeatletter
\def\eqref#1{(\ref{#1})}

\newcommand{\rright}{\right}
\newcommand{\lleft}{\left}
\newcommand{\rrVert}{\Vert}
\newcommand{\rrvert}{\vert}
\newcommand{\llVert}{\Vert}
\newcommand{\llvert}{\vert}

\newcommand{\sgn}{\@ifstar\sgnstar\sgnnostar}
\def\sgnnostar{\operatorname{sgn}}
\def\sgnstar{\overline{\operatorname{sgn}}}

\newtheorem{theorem}{Theorem}[section]
\newtheorem{proposition}{Proposition}[section]
\newremark{rem}{Remark}[section]
\makeatother

\begin{document}
\begin{frontmatter}

\title{Diffusions with rank-based characteristics and values in
the nonnegative quadrant}
\runtitle{Diffusions with rank-based characteristics and values in
the nonnegative quadrant}

\begin{aug}
\author[1]{\fnms{Tomoyuki} \snm{Ichiba}\thanksref{1}\ead[label=e1]{ichiba@pstat.ucsb.edu}},
\author[2]{\fnms{Ioannis} \snm{Karatzas}\thanksref{2}\ead[label=e2]{ik@enhanced.com}\ead[label=e4]{ik@math.columbia.edu}} \and
\author[3]{\fnms{Vilmos} \snm{Prokaj}\corref{}\thanksref{3}\ead[label=e3]{prokaj@cs.elte.hu}}
\runauthor{T. Ichiba, I. Karatzas and V. Prokaj} 
\address[1]{Department of Statistics and Applied Probability, South Hall,
University of California, Santa Barbara, CA~93106, USA. \printead{e1}}
\address[2]{INTECH Investment Management LLC, One Palmer Square,
Suite 441, Princeton, NJ 08542, USA and
Department of Mathematics, Columbia University, MailCode 4438, New
York, NY
10027, USA. \mbox{\printead{e2,e4}}}
\address[3]{Department of Probability Theory and Statistics, E\"otv\"os
Lor\'and University, 1117 Budapest, P\'azm\'any P\'eter s\'et\'any 1/C,
Hungary. \printead{e3}}
\end{aug}

\received{\smonth{2} \syear{2012}}
\revised{\smonth{6} \syear{2012}}

%
\begin{abstract}
We construct diffusions with values in the nonnegative orthant,
normal reflection along each of the axes, and two pairs of local drift/variance
characteristics assigned according to rank; one of the variances is
allowed to vanish,
but not both. The construction involves solving a system of coupled Skorokhod
reflection equations, then ``unfolding'' the Skorokhod reflection of a suitable
semimartingale in the manner of Prokaj (\textit{Statist. Probab.
Lett.} \textbf{79} (2009) 534--536). Questions of pathwise uniqueness and
strength are also addressed, for systems of stochastic differential
equations with reflection that realize these diffusions.
When the variance of the laggard is at least as large as that of the
leader, it is shown that the corner of the quadrant is never visited.
\end{abstract}

%

\end{frontmatter}

\section{Introduction}
\label{sec1}

We construct a planar diffusion $  \mathcal{X} (\cdot)= (X_1 (\cdot
), X_2 (\cdot))' $ according to the following recipe: each of the
components or ``particles'', $  X_1 (\cdot) $ and $  X_2 (\cdot)
$, starts at a nonnegative position, respectively, $  x_1\ge0 $ and
$  x_2 \ge0 $, and behaves locally like Brownian motion. The
characteristics of these motions are assigned not by name, but by rank:
the leader is assigned drift $  -h \le0 $ and dispersion $  \rho
\ge0 $, whereas the laggard is assigned drift $  g \ge0 $ and
dispersion $  \sigma> 0 $.
We force the planar process $  \mathcal{X} (\cdot)= ( X_1 (\cdot),
X_2 (\cdot))' $ never to leave the nonnegative quadrant in the
Euclidean plane, by imposing a reflecting barrier at the origin for the
laggard; this corresponds to orthogonal reflection along each of the
faces of the quadrant. In the interest of concreteness and simplicity,
we shall set
\begin{equation}
\label{1.1} \lambda :=  g + h  ,\qquad   \xi :=  x_1 +
x_2 >0 , \qquad   \rho^2 + \sigma^2  = 1  .
\end{equation}

A bit more precisely, we shall try to construct a filtered probability space
$  ( \Omega, \mathfrak{F}, \mathbb{P}),   \mathbf{F} = \{
\mathfrak{F} (t) \}_{0 \le t
< \infty} $ and on it two pairs $  (B_1 (\cdot), B_2 (\cdot)) $
and $
(X_1 (\cdot), X_2 (\cdot)) $ of continuous, $  \mathbf{F}$-adapted
processes, such that $  (B_1 (\cdot), B_2 (\cdot)) $ is planar Brownian
motion and $  (X_1 (\cdot), X_2 (\cdot)) $ a continuous
semimartingale that
takes values in the quadrant $  [0, \infty)^2 $ and satisfies the
dynamics
\begin{eqnarray}
\label{1.3} 
 \mathrm{d} X_1 (t)
 &= & ( g   \mathbf{ 1}_{ \{ X_1 (t) \le
X_2 (t)
\} } - h   \mathbf{ 1}_{ \{ X_1 (t) > X_2 (t) \} } )
\,\mathrm{d} t
\nonumber\\[-8pt]\\[-8pt]
&&{} +   \rho  \mathbf{ 1}_{ \{ X_1 (t) > X_2 (t) \} }  \,\mathrm{d} B_1 (t) +
\mathbf{ 1}_{ \{ X_1 (t) \le X_2 (t) \} }  \bigl( \sigma \, \mathrm{d} B_1 (t) +
\mathrm{d} L^{X_1 \wedge X_2} (t) \bigr) , \nonumber
\\
\label{1.4} 
 \mathrm{d} X_2 (t)
 &= & ( g   \mathbf{ 1}_{ \{ X_1 (t) > X_2 (t)
\} } - h   \mathbf{ 1}_{ \{ X_1 (t) \le X_2 (t) \} } )
\,\mathrm{d} t
\nonumber\\[-8pt]\\[-8pt]
&&{}+   \rho  \mathbf{ 1}_{ \{ X_1 (t) \le X_2 (t)
\} }  \,\mathrm{d} B_2
(t)
+  \mathbf{ 1}_{ \{ X_1 (t) > X_2 (t) \} }  \bigl( \sigma  \,\mathrm {d}
B_2 (t) + \mathrm{d} L^{X_1 \wedge X_2} (t) \bigr) .\nonumber
\end{eqnarray}
Here and in the sequel, we denote by $  L^{X } (\cdot) $ the local
time accumulated at the origin by a generic continuous semimartingale
$  X(\cdot) $, and by $  (X_1 \wedge X_2) (\cdot) = \min( X_1
(\cdot), X_2 (\cdot)) $, $  (X_1 \vee X_2) (\cdot) = \max( X_1
(\cdot), X_2 (\cdot)) $ the laggard and the leader, respectively, of
two such semimartingales $ X_1 (\cdot) $, $  X_2 (\cdot)  $. The
local time or ``boundary'' process $  L^{X_1 \wedge X_2} (\cdot) $
in (\ref{1.3}), (\ref{1.4}) imposes the reflecting boundary condition
on the laggard that we referred to earlier, and keeps the planar
process $  \mathcal{X} (\cdot)= (X_1 (\cdot), X_2 (\cdot))' $
from exiting the nonnegative quadrant. Because we are allowing one of
the two variances to be equal to zero, the system of equations (\ref
{1.3}), (\ref{1.4}) incorporates features of discontinuity,
degeneracy, and reflection on a nonsmooth boundary, all at once; this
makes its analysis challenging.

On a suitable filtered probability space, we shall construct fairly
explicitly a process $  \mathcal{X} (\cdot)= ( X_1 (\cdot),$ $ X_2
(\cdot))'  $ with continuous paths and values in the quadrant $  [0,
\infty)^2$, along with a planar Brownian motion $ (B_1 (\cdot), B_2
(\cdot)) $, so that the equations (\ref{1.3}), (\ref{1.4}) and the
following properties are satisfied $  \mathbb{P}$-a.e., the last one
for $
(i, j) \in\{ (1,2), (2,1) \} $:
\begin{eqnarray}
 \label{1.6}
 \int_0^\infty\mathbf{
1}_{\{ X_1 (t) = X_2 (t)\}} \, \mathrm{d} t   &=& 0 , \qquad   \int_0^\infty
\mathbf{ 1}_{\{ X_1 (t) \wedge X_2 (t) =0\}}   \,\mathrm {d} t  = 0 ,
\\
\label{1.7}
\int_0^\infty\mathbf{
1}_{\{ X_1 (t) = X_2 (t)\}}  \,\mathrm{d} L^{X_1 \wedge
X_2} (t)  &=& 0 , \qquad    \int
_0^\infty\mathbf{ 1}_{\{ X_1 (t) \wedge X_2 (t)= 0\}}  \,\mathrm {d}
L^{X_i - X_j} (t)  = 0  .
\end{eqnarray}
In a terminology first introduced apparently by Manabe and Shiga
\cite{MR0346907}, the second condition in (\ref{1.6}) mandates that
the faces of the quadrant are ``nonsticky" for the planar diffusion $
\mathcal{X} (\cdot)= (X_1 (\cdot), X_2 (\cdot))' $; whereas the first
condition in (\ref{1.6}) can be interpreted as saying that the
diagonal of
the quadrant is also ``nonsticky'' for this diffusion.

The conditions of (\ref{1.7}) can be interpreted, in the spirit of Reiman and Williams
\cite{MR921820}, as saying that the ``boundary
processes'' do not charge the set of times when the diffusion is at the
intersection of two faces. We show in Proposition \ref{Prop0} that the
properties of (\ref{1.7}) are satisfied automatically, as long as the
process $  \mathcal{X} (\cdot) $ stays away from the corner of the
quadrant.

We shall prove the following results, Theorems \ref{Thm1}--\ref{Thm2}
below. In Theorem \ref{Thm2}, we shall impose
\begin{equation}
\label{SigmaGeRho} 1/2   \le  \sigma^2   \le  1 ,
\end{equation}
a condition mandating that the variance of the laggard be at least as
big as that of the leader. Under this condition, it will turn out that
the two particles never collide with each other at the origin, so the
process $  \mathcal{X} (\cdot) $ takes values in the punctured
nonnegative quadrant
\begin{equation}
\label{1.2} \mathfrak{ S}  :=  [ 0, \infty)^2 \setminus\bigl\{
(0,0)\bigr\} .
\end{equation}

\begin{theorem}
\label{Thm1}
The system of stochastic differential equations (\ref{1.3}) and (\ref
{1.4}) admits a weak solution, with a state-process $  \mathcal{X}
(\cdot)= ( X_1 (\cdot), X_2 (\cdot))' $ that takes values in the
quadrant $  [0, \infty)^2 $ and satisfies the properties of (\ref
{1.6}), (\ref{1.7}); and when restricted to the stochastic interval
$  [0, \tau) $, with
\begin{equation}
\label{1.5} \tau :=  \inf\bigl\{   t > 0 \dvt  X_1 (t) =
X_2 (t) =0  \bigr\}
\end{equation}
the first hitting time of the corner of the quadrant, this solution is
pathwise unique, thus also strong.
\end{theorem}

\begin{theorem}
\label{Thm3}
In the nondegenerate case $  \sigma<1  $, and among weak solutions
that satisfy the conditions of (\ref{1.7}), the solution of Theorem
\ref{Thm1} is unique in distribution with a continuous, strongly
Markovian state-process $  \mathcal{X} (\cdot) $.
\end{theorem}

\begin{theorem}
\label{Thm2}
Under the condition (\ref{SigmaGeRho}), the hitting time of (\ref
{1.5}) is a.s. infinite: $  \mathbb{P}(\tau= \infty) =1.$ The
system of
(\ref{1.3}) and (\ref{1.4}) admits then a pathwise unique, strong solution.
\end{theorem}

\subsection{Preview}
\label{ }

A weak solution of the system (\ref{1.3})--(\ref{1.4}) is constructed rather
explicitly in Section \ref{sec4} (Proposition \ref{Prop2}) following an
\textit{a priori} analysis of its structure in Section \ref{sec3}, and is shown
to satisfy the properties (\ref{1.6}) and (\ref{1.7}). This construction
involves finding (in Section \ref{sec5}, proof of Proposition~\ref{Prop1})
the unique solution of a system of coupled \textsc{Skorokhod} reflection
equations; and the ``unfolding'', in Section \ref{sec4.1.a}, of an
appropriate nonnegative semimartingale in the manner of
Prokaj \cite{MR2494646}. The constructed solution admits the \textit{skew
representations} of (\ref{2.23b})--(\ref{2.24}).

It is shown that the constructed state process $  \mathcal{X} (\cdot
) $ never visits the corner of the quadrant under the condition (\ref
{SigmaGeRho}) (Section \ref{sec6}, Proposition \ref{Prop3}). When
this condition fails, the process $  \mathcal{X} (\cdot) $ can hit
the corner of the quadrant; but then ``it knows how to extricate
itself'' in such a manner that uniqueness in distribution holds, as
shown in Section \ref{subsec8} (proof of Theorem \ref{Thm3}).

Pathwise uniqueness is established in Section \ref{sec7}, Proposition
\ref{Prop4}. Questions of pathwise uniqueness and strength, for
additional systems of stochastic differential equations with reflection
that realize this diffusion, are addressed in Section \ref{sec9}.
Issues of recurrence, transience and invariant densities are touched
upon briefly in Section \ref{sec10}. The degenerate case $  \sigma
=0 $ (zero variance for the laggard) is discussed briefly in the
\hyperref[app]{Appendix}. Basic facts about semimartingale local
time are recalled in Section \ref{sec2a}.

\subsection{Connections}
\label{ }

Diffusions with rank-based characteristics were introduced in Fernholz \cite{MR1894767}
and studied in Banner, Fernholz and Karatzas \cite{MR2187296}, Ichiba \textit{et al.} \cite{2009arXiv0909.0065I}
in connection with the study of long-term stability properties of large equity
markets. Their detailed probabilistic study includes
Fernholz \textit{et al.} \cite{2011arXiv1108.3992F} in two dimensions, and Ichiba \textit{et~al.} \cite{2011arXiv1109.3823I}
in three or more dimensions. This paper extends the results of
Fernholz \textit{et al.} \cite{2011arXiv1108.3992F} to a situation where -- through reflection
at the
faces of the nonnegative quadrant -- the vector of ranked processes $  (X_1
(\cdot) \vee X_2 (\cdot) , X_1 (\cdot) \wedge X_2 (\cdot)) $ has
itself a
stable distribution, under the conditions $  h>g>0 $; cf. Section
\ref{sec10.1} for some explicit computations. This has important ramifications
for parameter estimation via time-reversal, as explained in
Fernholz, Ichiba and Karatzas \cite{springerlink:10.1007/s10436-012-0193-2}.
It is an interesting question, whether the analysis in this paper can
be extended, to study multidimensional diffusions with rank-based
characteristics and reflection on the faces of an orthant.

There are also rather obvious connections of the model studied here with
queueing models of the so-called ``generalized Jackson type'' under
heavy-traffic conditions; we refer the reader to
Foschini \cite{MR651627}, Harrison and Williams \cite{MR912049}, and Reiman \cite{MR757317}.

\section{On semimartingale local time}
\label{sec2a}


Let us start with a continuous, real-valued semimartingale
\begin{equation}
\label{Semi} X (\cdot)  =  X (0) + \Theta (\cdot) +C(\cdot)  ,
\end{equation}
where $  \Theta (\cdot) $ is a continuous local martingale and $  C
(\cdot) $
a continuous process of finite first variation such that $  \Theta
(0) =
C(0) =0 $; note that $ \langle X \rangle(\cdot)=\langle\Theta
\rangle(\cdot) $. The \textit{local time} $  L^{ X } (t) $ accumulated at the origin over the time-interval $  [0,t] $ by this process, is given as
\begin{equation}
\label{Tanaka} L^{ X} (t) :=  \lim_{\varepsilon\downarrow0}   { 1 \over 2
  \varepsilon } \int_0^t
\mathbf{1}_{ \{ 0 \le X (s) <
\varepsilon\} }  \, \mathrm{d} \langle X \rangle 
( s) =
\bigl( X (t) \bigr)^+ - \bigl( X (0) \bigr)^+ - \int_0^t
\mathbf{ 1}_{\{ X (s) > 0 \} } \, \mathrm{d} X (s)  .
\end{equation}
This defines a nondecreasing, continuous and adapted process $  L^{ X}
( t) $, $ 0 \le t < \infty $ which is flat off the zero set of $  X
(\cdot)  $, namely
\begin{eqnarray}
\label{flat} \int_0^\infty  \mathbf{
1}_{ \{ X (t) \neq0\} } \, \mathrm{d} L^{
X} ( t)  &=& 0    ; \qquad  \mbox{we
also have the property}
\nonumber\\[-8pt]\\[-8pt]
    \int_0^\infty
\mathbf{ 1}_{ \{ X (t) = 0\} } \, \mathrm{d} \langle X \rangle(t)  &=&  0 .\nonumber
\end{eqnarray}

\noindent$\bullet $
On the other hand, for a \textit{nonnegative} continuous semimartingale
$  X (\cdot)  $ of the form (\ref{Semi}), we obtain from (\ref
{Tanaka}), (\ref{flat}) the representations
\begin{equation}
\label{2.19.a} L^{ X} (\cdot) =  \int_0^\cdot
  \mathbf{ 1}_{ \{ X (t) = 0\} }\,  \mathrm{d} X (t)  =  \int_0^\cdot
  \mathbf{ 1}_{ \{ X (t) = 0\}
} \, \mathrm{d} C (t) .
\end{equation}
Using this observation, it can be shown as in
Ouknine \cite{MR1071565} (see also Ouknine and Rutkowski \cite{MR1324099})
that the local time at the origin of the laggard of two continuous,
nonnegative semimartingales $  X_1(\cdot) $, $  X_2(\cdot) $ is
given as
\begin{eqnarray}
\label{Ouk4} %
  L^{X_1 \wedge X_2} ( t)  &=&  \int
_0^t \mathbf{ 1}_{\{ X_1 (s) \le
X_2 (s)
\}}  \, \mathrm{d}
L^{X_1 } (s)  + \int_0^t \mathbf{
1}_{\{ X_1 (s) > X_2 (s)
\}}  \, \mathrm{d} L^{X_2 } (s)
\nonumber\\[-8pt]\\[-8pt]
&&{} - \int_0^t \mathbf{ 1}_{\{ X_1 (s) = X_2 (s)=0 \}}
\,\mathrm{d} L^{X_1 - X_2} (s);\nonumber    %
\end{eqnarray}
a companion representation, in which the r\^oles of $ X_1 $ and $
X_2 $
are interchanged, is also valid. With the help of (\ref{Ouk4}),
Ouknine \cite{MR1071565} derives a purely algebraic proof of the
Yan \cite{MR619601,MR834059} identity
\begin{equation}
\label{Ouk} L^{X_1 \vee X_2} ( \cdot) + L^{X_1 \wedge X_2} ( \cdot)  =
L^{X_1
} ( \cdot) + L^{ X_2} ( \cdot) .
\end{equation}

\subsection{Tanaka formulae}
\label{sec2aaa}

For a continuous, real-valued semimartingale $  X(\cdot) $ as in
(\ref{Semi}), and with the conventions
\begin{equation}
\label{sgn} \overline{\operatorname{sgn}}  (x)  :=  \mathbf{ 1}_{(0, \infty)} (x)
- \mathbf{ 1}_{(- \infty, 0)} (x)  , \qquad   \operatorname{sgn}  (x)  :=  \mathbf{
1}_{(0, \infty)} (x) - \mathbf{ 1}_{(- \infty, 0]} (x) , \qquad   x \in
\mathbb{R}
\end{equation}
for the symmetric and the left-continuous versions of the signum
function, the \textsc{Tanaka} formula
\begin{equation}
\label{Tanaka2} \big| X (\cdot)\big|  =  \big| X (0)\big| + \int_0^\cdot
\operatorname{sgn} \bigl( X(t) \bigr)  \, \mathrm{d} X(t)  +  2  L^X (
\cdot)
\end{equation}
holds. Applying this formula to the continuous, nonnegative
semimartingale $  | X(\cdot)| $, then comparing with the expression
of (\ref{Tanaka2}) itself, we obtain the generalization
\begin{equation}
\label{2.19.c} 2  L^{ X} (\cdot) - L^{ |X|} (\cdot) =  \int
_0^\cdot  \mathbf{ 1}_{ \{ X (t) = 0\} }
\,\mathrm{d} X (t)  =  \int_0^\cdot  \mathbf{
1}_{ \{ X (t) = 0\} } \, \mathrm{d} C (t)
\end{equation}
of the representation (\ref{2.19.a}), and from it the companion
\textsc{Tanaka} formula
\begin{equation}
\label{Tanaka1} \big| X (\cdot)\big|  =  \big| X (0)\big| + \int_0^\cdot
\overline{\operatorname{sgn}} \bigl( X(t) \bigr)  \, \mathrm{d} X(t)  +
L^{|X|} (\cdot) .
\end{equation}
It follows from (\ref{2.19.c}) that the identity $  L^{ |X|} (\cdot)
\equiv
2  L^{ X} (\cdot)  $ holds when $  \int_0^\cdot  \mathbf{ 1}_{
\{ X
(t) = 0\} }  \,\mathrm{d} C (t) \equiv0 $. (This last condition
guarantees also the continuity of the ``local time random field'' $  a
\longmapsto L^X (\cdot , a) $ in its spatial argument at the origin
$
a^*=0 $; cf. page 223 in
Karatzas and Shreve \cite{MR1121940}.)\vadjust{\goodbreak} In light of (\ref{flat}), the identity $  L^{ |X|}
(\cdot) \equiv2  L^{ X} (\cdot)  $ holds when the finite variation
process $  C(\cdot) $ is absolutely continuous with respect to the
bracket $
\langle X \rangle 
(\cdot) $ of the local martingale part of the semimartingale.

For the theory that undergirds these results we refer, for instance, to
Karatzas and Shreve \cite{MR1121940}, Section 3.7.

\subsection{Ramifications}
\label{sec2aa}

For two continuous, \textit{nonnegative} semimartingales $  X_1(\cdot),
 X_2 (\cdot)  $ that satisfy the last properties in (\ref{1.7}),
the expression of (\ref{Ouk4}) gives the $  \mathbb{P}$-a.e. representation
\begin{equation}
\label{Ouk5} L^{X_1 \wedge X_2} ( t)  =  \int_0^t
\mathbf{ 1}_{\{ X_1 (s) \le
X_2 (s) \}}  \, \mathrm{d} L^{X_1 } (s) + \int
_0^t \mathbf{ 1}_{\{
X_1 (s) > X_2 (s) \}}  \, \mathrm{d}
L^{X_2 } (s)
\end{equation}
for the local time of the laggard. On the strength of the properties
(\ref{flat}) and (\ref{1.7}), we deduce then
\begin{eqnarray}
\label{Ouk1} \int_0^t \mathbf{
1}_{\{ X_1 (s) \le X_2 (s) \}} \,  \mathrm{d} L^{X_1 \wedge X_2} ( s)  &= & \int
_0^t \mathbf{ 1}_{\{ X_1 (s) \le
X_2 (s) \}}   \,\mathrm{d}
L^{X_1 } (s)  =  L^{X_1} (t) ,
\nonumber\\[-8pt]\\[-8pt]
\int_0^t \mathbf{
1}_{\{ X_1 (s) = X_2 (s) \}}  \, \mathrm{d} L^{X_1
} (s) &=&  \int
_0^t \mathbf{ 1}_{\{ X_1 (s) = X_2 (s) \}}  \, \mathrm {d}
L^{X_1 \wedge X_2} ( s)  = 0  ,\nonumber
\end{eqnarray}
thus also $  \int_0^t \mathbf{ 1}_{\{ X_1 (s) = X_2 (s) \}}
\,\mathrm{d} L^{X_2} (s) =0 $
by symmetry and
\begin{eqnarray}
\label{Ouk2} %
  \int_0^t
\mathbf{ 1}_{\{ X_1 (s) > X_2 (s) \}}  \, \mathrm{d} L^{X_1
\wedge
X_2} ( s)  &=&  \int
_0^t \mathbf{ 1}_{\{ X_1 (s) > X_2 (s) \}} \,  \mathrm{d}
L^{X_2 } (s)
\nonumber\\[-8pt]\\[-8pt]
 &=&  \int_0^t \mathbf{ 1}_{\{ X_1 (s) \ge X_2 (s) \}}
 \, \mathrm {d} L^{X_2 } (s)  =  L^{X_2} (t)\nonumber
\end{eqnarray}
for $  0 \le t < \infty $. It follows from (\ref{Ouk1}), (\ref
{Ouk2}) and (\ref{Ouk}) that we have then
\begin{equation}
\label{Ouk3} L^{X_1 \wedge X_2} ( \cdot)  =  L^{X_1 } ( \cdot) +
L^{ X_2} ( \cdot) ,\qquad  \mbox{whence }   L^{X_1 \vee X_2} ( \cdot)
 =  0 .
\end{equation}
To wit: for any two continuous, nonnegative semimartingales $  X_1
(\cdot) $ and $  X_2 (\cdot)  $ that satisfy the properties of
(\ref{1.7}), the leader $  X_1 (\cdot) \vee X_2 (\cdot)  $ does
not accumulate any local time at the origin, \textit{even in situations}
(such as in Proposition \ref{Prop5} below) \textit{where the planar
process $  (X_1 (\cdot), X_2 (\cdot)) $ does attain the corner of
the quadrant}.

\noindent$\bullet $
It is fairly clear from this discussion, in particular from (\ref
{Ouk1}) and (\ref{Ouk2}), that, in the presence of (\ref{1.7}), the
dynamics of (\ref{1.3}), (\ref{1.4}) can be cast in the more
``conventional'' form
\begin{eqnarray}
\label{1.3.a} %
  \mathrm{d} X_1 (t) &=  &
({g   \mathbf{ 1}_{ \{ X_1 (t) \le
X_2 (t) \} } - h   \mathbf{ 1}_{ \{ X_1 (t) > X_2 (t) \} }} )
\,\mathrm{d} t
\nonumber\\[-8pt]\\[-8pt]
&&{} +   \rho  \mathbf{ 1}_{ \{ X_1 (t) > X_2 (t) \} } \, \mathrm{d} B_1 (t) +
\sigma  \mathbf{ 1}_{ \{ X_1 (t) \le X_2 (t) \} }  \,\mathrm{d} B_1 (t) +
\mathrm{d} L^{X_1 } (t)  , \nonumber  %
\\
\label{1.4.a} %
  \mathrm{d} X_2 (t) &= & ({g
  \mathbf{ 1}_{ \{ X_1 (t) > X_2
(t) \} } - h   \mathbf{ 1}_{ \{ X_1 (t) \le X_2 (t) \} } } )
\,\mathrm{d} t
\nonumber\\[-8pt]\\[-8pt]
&&{} +   \rho  \mathbf{ 1}_{ \{ X_1 (t) \le X_2 (t) \} }  \,\mathrm{d} B_2 (t) +
\sigma  \mathbf{ 1}_{ \{ X_1 (t) > X_2 (t) \} } \, \mathrm{d} B_2(t) + \mathrm{d}
L^{ X_2} (t)  .\nonumber   %
\end{eqnarray}
Here each of the components $  X_1 (\cdot) $, $  X_2 (\cdot) $ of
the planar process $  \mathcal{X} (\cdot) $ is reflected at the
origin via its own local time, respectively $  L^{X_1} (\cdot)  $
and $  L^{X_2} (\cdot)  $. Conversely, in the presence of condition
(\ref{1.7}), the dynamics of (\ref{1.3.a})--(\ref{1.4.a}) lead to
those of (\ref{1.3})--(\ref{1.4}).

\begin{proposition}
\label{Prop0}
For a continuous, planar semimartingale $  \mathcal{X}(\cdot) =
(X_1(\cdot), X_2 (\cdot))'  $ that takes values in the punctured
quadrant $ \mathfrak{ S}  $ of (\ref{1.2}), the conditions of (\ref
{1.7}) are satisfied automatically.
\end{proposition}

\begin{pf}
In this case, and in conjunction with (\ref{flat}), (\ref{2.19.a}),
the expression (\ref{Ouk4}) takes the form
\begin{equation}
\label{1.7.a} L^{X_1 \wedge X_2} ( t)  =  \int_0^t
\mathbf{ 1}_{\{ X_2 (s) > X_1
(s)=0 \}}  \, \mathrm{d} L^{X_1 } (s) + \int
_0^t \mathbf{ 1}_{\{ X_1
(s) > X_2 (s) =0\}} \,  \mathrm{d}
L^{X_2 } (s)  ,
\end{equation}
and
the first property in (\ref{1.7}) follows; the second is then a
consequence of (\ref{flat}), (\ref{2.19.a}).
\end{pf}

\section{Analysis}
\label{sec3}

Let us suppose that such a probability space as stipulated in Section
\ref{sec1} has been constructed, and on it a pair $  B_1 (\cdot) $,
$ B_2
(\cdot) $ of independent standard Brownian motions, as well as two
continuous, nonnegative semimartingales $  X_1(\cdot),  X_2 (\cdot)
 $
such that \textit{the dynamics \textup{(\ref{1.3})--(\ref{1.4})} and the
conditions of
\textup{(\ref{1.7})} are satisfied.} We import the notation of
Fernholz \textit{et al.} \cite{2011arXiv1108.3992F}: in addition to (\ref{1.1}), we set
\begin{equation}
\label{2.1} \nu =  g-h , \qquad   y  =  x_1 - x_2
 ,  \qquad  r_1 =  x_1 \vee x_2 ,\qquad
r_2 =  x_1 \wedge x_2 ,
\end{equation}
and introduce the difference and the sum of the two component
processes, namely
\begin{equation}
\label{2.2} Y (\cdot)  :=  X_1 (\cdot) - X_2 (\cdot)
 ,\qquad   \Xi(\cdot)   :=  X_1 (\cdot) + X_2 (\cdot) .
\end{equation}
We introduce also the two planar Brownian motions $   ( W_1 (\cdot
), W_2 (\cdot)  ) $ and $   ( V_1 (\cdot), V_2 (\cdot)
) $, given by
\begin{eqnarray}
\label{2.3} W_1 (t)  &:=&  \int_0^t
\mathbf{ 1}_{ \{ Y (s) > 0 \} }  \, \mathrm {d} B_1 (s) - \int
_0^t \mathbf{ 1}_{ \{ Y (s) \le0 \} } \,  \mathrm{d}
B_2 (s)   ,
\\
\label{2.4} W_2 (t)  &:= & \int_0^t
\mathbf{ 1}_{ \{ Y (s) \le0 \} } \,  \mathrm {d} B_1 (s) - \int
_0^t \mathbf{ 1}_{ \{ Y (s) > 0 \} } \,  \mathrm{d}
B_2 (s)
\end{eqnarray}
and
\begin{eqnarray}
\label{2.5} V_1(t)  &:=&  \int_0^t
\mathbf{ 1}_{ \{ Y (s) > 0 \} }  \, \mathrm {d} B_1 (s) + \int
_0^t \mathbf{ 1}_{ \{ Y (s) \le0 \} } \,  \mathrm{d}
B_2 (s)  ,
\\
\label{2.6} V_2(t)  &:=&   \int_0^t
\mathbf{ 1}_{ \{ Y (s) \le0 \} } \,  \mathrm {d} B_1 (s) + \int
_0^t \mathbf{ 1}_{ \{ Y (s) > 0 \} }  \, \mathrm{d}
B_2 (s)  ,
\end{eqnarray}
respectively. Finally, we construct the Brownian motions $  W (\cdot)
$, $  V (\cdot)  $, $  Q (\cdot) $ and $  V^\flat(\cdot)   $ as
\begin{eqnarray}
 \label{2.7} W (\cdot) &:= &\rho  W_1 (\cdot) + \sigma
W_2 (\cdot) , \qquad   V (\cdot) := \rho  V_1 (\cdot)
+ \sigma  V_2 (\cdot) ,
\\
\label{2.9.b} Q (\cdot) &:=&\sigma  V_1 (\cdot) + \rho
V_2 (\cdot)  , \qquad   V^{ \mathbf{ \flat}}( \cdot) := \rho
V_1 (\cdot) - \sigma  V_2 (\cdot) ,
\end{eqnarray}
note that $  Q (\cdot) $ and $  V^\flat(\cdot)   $ are
independent, and observe the intertwinements
\begin{eqnarray}
\label{2.10} V_1(t) &=& \int_0^t
\operatorname{sgn} \bigl( Y (s) \bigr) \,\mathrm{d} W_1 (s) ,\qquad
V_2(t) = -\int_0^t \operatorname{sgn}
\bigl( Y (s) \bigr) \,\mathrm{d} W_2 (s)  ,
\nonumber\\[-8pt]\\[-8pt]
   V^\flat(t) &=&
\int_0^t \operatorname{sgn} \bigl( Y (s) \bigr)
\,\mathrm{d} W (s) .\nonumber
\end{eqnarray}

\noindent$\bullet $ After this preparation, and recalling the first property of
(\ref{1.7}), we observe that the difference $  Y(\cdot) $ and the
sum $  \Xi(\cdot) $ from (\ref{2.2}) satisfy, respectively, the equations
\begin{eqnarray}\label{2.12}
\label{2.11} Y (t)  &=&  y + \int_0^t
\operatorname{sgn} \bigl( Y(s) \bigr)  \bigl( - \lambda \, \mathrm{d} s - \mathrm{d}
L^{X_1 \wedge X_2}(s)+ \mathrm{d}  V^{\flat}(s) \bigr)
\nonumber\\[-8pt]\\[-8pt]
&=&  y + \int_0^t \operatorname{sgn}
\bigl( Y(s) \bigr)  \bigl( - \lambda \,  \mathrm{d} s - \mathrm{d}
L^{X_1 \wedge X_2}(s) \bigr) + W(t)\nonumber
\end{eqnarray}
and
\begin{equation}
  \Xi(t)  = \xi+
\nu t + V (t) + L^{X_1 \wedge X_2}(t) , \qquad  0 \le t < \infty
\end{equation}
in the notation of (\ref{1.1}), (\ref{2.1}). An application of the
\textsc{Tanaka} formula (\ref{Tanaka2}) to the semimartingale $
Y(\cdot)$ of (\ref{2.11}) represents now the size of the ``gap''
between $  X_1 (t) $ and $  X_2 (t) $ as
\begin{equation}
\label{2.13} \big| Y (t) \big|  =| y | - \lambda t - L^{X_1 \wedge X_2}(t) +
V^\flat(t) + 2  L^Y (t) ,\qquad   0 \le t < \infty .
\end{equation}
On the other hand, with the help of the theory of the \textsc{Skorokhod}
reflection problem (e.g.,
Karatzas and Shreve \cite{MR1121940}, page 210), we obtain from (\ref{2.13}), (\ref
{flat}) the equation
\begin{equation}
\label{2.14} 2   L^Y (t)  =  \max_{0 \le s \le t} \bigl( - |y| +
\lambda  s + L^{X_1 \wedge X_2}(s) - V^{ \flat} (s) \bigr)^+ , \qquad  0 \le t
< \infty .
\end{equation}

\subsection{Ranks}
\label{sec3.1}

It is convenient now to introduce explicitly the ranked versions
\begin{equation}
\label{ranks} R_1 (\cdot)   :=   X_1 (\cdot) \vee
X_2 (\cdot)   , \qquad   R_2 (\cdot) :=
X_1 (\cdot) \wedge X_2 (\cdot)
\end{equation}
of the components of the vector process $  \mathcal{X} (\cdot)= (X_1
(\cdot), X_2 (\cdot))'  $. From (\ref{2.12}), (\ref{2.13}), we have
\begin{eqnarray}\label{2.15.a}
R_1 (t) + R_2 (t) &=&  X_1 (t)
+ X_2 (t)  =  \xi+ \nu t + V(t) + L^{R_2} (t) ,\qquad   0
\le t < \infty,
\nonumber\\[-8pt]\\[-8pt]
R_1 (t) - R_2 (t) &=&   \big|
X_1 (t) - X_2 (t) \big| =   \big| Y (t) \big|  =  |y| +
V^{ \flat} (t) - \lambda  t - L^{R_2}(t) + 2  L^Y
(t) ,\quad \nonumber
\end{eqnarray}
and these representations lead to the expressions
\begin{eqnarray}
\label{2.16} R_1 (t) &=&  r_1 - h  t + \rho
V_1(t) + L^Y (t) , \qquad 0 \le t < \infty,
\\
\label{2.17}   R_2 (t) &=&  r_2 + g  t + \sigma
V_2(t)- L^Y (t) + L^{R_2} (t) , \qquad 0 \le t <
\infty .
\end{eqnarray}

A few remarks are in order. The equations (\ref{2.16}), (\ref{2.17})
identify the processes $  V_1(\cdot) $ and $  V_2(\cdot) $ of
(\ref{2.5}), (\ref{2.6}) as the independent Brownian motions
associated with individual \textit{ranks}, the ``leader'' $  R_1 (\cdot
) $ and the ``laggard'' $  R_2 (\cdot) $, respectively; whereas the
independent Brownian motions $  B_1(\cdot) $ in (\ref{1.3}) and $
B_2(\cdot) $ in (\ref{1.4}) are associated with the specific
``names'' (indices, or identities) of the individual particles. On the
other hand the equation (\ref{2.17}) leads, in conjunction with (\ref
{flat}) and the theory of the \textsc{Skorokhod} reflection problem
once again, to the representation
\begin{equation}
\label{2.18} L^{R_2} (t)  =  \max_{0 \le s \le t} \bigl( -
r_2 - g   s + L^{Y}(s) - \sigma  V_{2} (s)
\bigr)^+ , \qquad  0 \le t < \infty .
\end{equation}

\noindent$\bullet $
Let us apply the second observation in (\ref{flat}) to the nonnegative
semimartingale $  R_1 (\cdot) - R_2 (\cdot) = | Y (\cdot)|  $ in
(\ref{2.15.a}); we obtain the first property of (\ref{1.6}), that is
\begin{equation}
\label{2.8} \int_0^\cdot  \mathbf{
1}_{ \{ X_1 (t) = X_2 (t)\} }\, \mathrm{d} \bigl\langle V^{\flat} \bigr\rangle(t)  =
\int_0^\cdot  \mathbf{ 1}_{ \{
X_1 (t) = X_2 (t)\} }
\,\mathrm{d} t  =  0  ,
\end{equation}
therefore also $  \int_0^\cdot  \mathbf{ 1}_{ \{ X_1 (t) = X_2
(t)\} }\, \mathrm{d} V^{\flat} (t)= 0 $. Whereas the observation
(\ref{2.19.a}) leads to
\begin{equation}
\label{2.19} L^{R_1-R_2} (\cdot) = \int_0^\cdot
  \mathbf{ 1}_{ \{ X_1 (t) = X_2
(t)\} }  \bigl( - \lambda \, \mathrm{d} t - \mathrm{d}
L^{R_2}(t) + \mathrm{d}  V^{\flat}(t) + 2 \,\mathrm{d}
L^{Y}(t) \bigr) = 2   L^Y (\cdot) .
\end{equation}
For this last identity, we have used the fact that $ L^Y (\cdot) $ is
supported on the set $ \{ t \ge0\dvt Y (t) = 0 \}$ $ = \{   t \ge0\dvt X_1
(t) = X_2 (t)  \} ,$ thanks to (\ref{flat}); that this set has zero
\textsc{Lebesgue} measure (the condition (\ref{2.8})); and that the
local time $  L^{R_2} (\cdot) $ is flat on this set (from the first
property in (\ref{1.7})).

Finally, we observe that the second property in (\ref{flat}), applied
to the nonnegative semimartingale $  R_2 (\cdot)  $ of (\ref
{2.17}), yields the second property in (\ref{1.6}).

\begin{rem}
\label{22}
In light of (\ref{2.19}), the equations (\ref{2.16}), (\ref{2.17})
assume the more suggestive form
\begin{eqnarray}
\label{2.20} R_1 (t) &=&  r_1 - h  t + \rho
V_1(t) + \tfrac{ 1}{2 } L^{R_1 - R_2} (t) , \qquad  0 \le t < \infty,
\\
\label{2.21}  R_2 (t) &=&  r_2 + g  t + \sigma
V_2(t)- \tfrac{ 1}{2 } L^{R_1 -
R_2} (t) + L^{R_2} (t)
 , \qquad  0 \le t < \infty .
\end{eqnarray}

In the nondegenerate case $  \sigma<1  $, these equations (\ref
{2.20}), (\ref{2.21}) give the filtration comparisons
\begin{equation}
\label{flit} \mathfrak{F}^{ (V_1, V_2)} (t)   \subseteq
\mathfrak{F}^{ (R_1, R_2)} (t)  ,\qquad   0 \le t < \infty .
\end{equation}
\end{rem}

\subsection{Skew representations}
\label{sec3.2}

On the strength of the representations (\ref{2.10}), (\ref{2.15.a})
and the notation of (\ref{1.1}) and (\ref{2.9.b}), the Brownian
motion $  V(\cdot) $ in (\ref{2.7}) can be cast in the form
\begin{eqnarray*}
V (t) &=&  \bigl(\rho^2 - \sigma^2 \bigr)
 V^\flat(t) + 2 \rho  \sigma Q (t)
\\
 &=&  \bigl(\rho^2-\sigma^2 \bigr)  \bigl(  \big|Y (t)
\big|-|y| + \lambda  t+L^{R_2}(t)-2 L^Y(t)  \bigr) + 2 \rho
 \sigma  Q (t) .
\end{eqnarray*}
In conjunction with the equations $ X_1 (t) + X_2 (t)  =  \xi+ \nu
 t + V(t) + L^{R_2} (t) $ and $   X_1 (t) -  X_2 (t)  = Y(t) $,
and with the notation
\begin{equation}
\label{2.22} \mu :=  \tfrac12 \bigl( \nu+ \lambda  \bigl( \rho^2 -
\sigma^2 \bigr) \bigr)  =  g  \rho^2 - h
\sigma^2  ,
\end{equation}
this leads \textit{for all} $  t \in[0, \infty) $ to the \textit{skew
representations}
\begin{eqnarray}
\label{2.23b} X_1 (t) &=& x_1 + \mu  t +
\rho^2 \bigl( Y^+ (t) - y^+ \bigr) - \sigma^2 \bigl( Y^-
(t) - y^- \bigr)
\nonumber\\[-8pt]\\[-8pt]
&&{} + \rho\sigma Q (t) + \bigl( \rho^2 -
\sigma^2 \bigr) \bigl( \tfrac{ 1}{2 }  L^{R_2} (t) -
L^Y (t) \bigr) ,\nonumber
\\
\label{2.24} X_2 (t) &=& x_2 + \mu  t -
\sigma^2   \bigl( Y^+ (t) - y^+ \bigr) + \rho^2 \bigl(
Y^- (t) - y^- \bigr)
\nonumber\\[-8pt]\\[-8pt]
&&{}  + \rho\sigma Q (t)+ \bigl( \rho^2 -
\sigma^2 \bigr) \bigl( \tfrac{ 1}{2 }  L^{R_2} (t) -
L^Y (t) \bigr) .\nonumber
\end{eqnarray}

\section{Synthesis}
\label{sec4}

Let us fix real constants $  x_1 \ge0 $, $  x_2 \ge0 $, $  h \ge
0$, $  g \ge0 $, $  \rho\ge0 $ and $ \sigma>0 $, and recall the
notation and assumptions of (\ref{1.1}), (\ref{2.1}).
We start with a filtered probability space $  ( \widetilde{\Omega},
\widetilde{\mathfrak{F}}, \widetilde{\mathbb{P}}),   \widetilde
{\mathbf
{F}} = \{ \widetilde{\mathfrak{F}} (t) \}_{0 \le t < \infty} $ and
two independent Brownian motions $  V_1 (\cdot), V_2 (\cdot)  $,
use these to create additional Brownian motions
\begin{eqnarray}
\label{3.1.a} V (\cdot)  &:=&   \rho  V_1 (\cdot) + \sigma
V_2 (\cdot) , \qquad  Q (\cdot)  :=  \sigma  V_1 (\cdot)+
\rho  V_2 (\cdot)  ,
\\
\label{3.1.b} V^{ \mathbf{ \flat}}( \cdot)  &:=&  \rho  V_1 (\cdot) -
\sigma   V_2 (\cdot)  , \qquad  Q^{ \mathbf{ \flat}}( \cdot)  :=  \sigma
V_1 (\cdot) - \rho   V_2 (\cdot) ,
\end{eqnarray}
as in (\ref{2.7}), (\ref{2.9.b}), and note that $  V (\cdot),
Q^{ \mathbf{ \flat}}( \cdot) $ are independent; the same is true
of $  Q (\cdot),  V^{ \mathbf{ \flat}}( \cdot) $.

With these ingredients we construct two continuous, increasing and
adapted processes $  A(\cdot) $, $  \Lambda(\cdot) $ with $
A(0) = \Lambda(0) =0 $ that satisfy for $  0 \le t < \infty $ the
system of equations
\begin{eqnarray}
\label{3.2.a} 2  A (t)  &=&  \max_{0 \le s \le t} \bigl(-|y|+\lambda s +
\Lambda(s)-V^{ \flat
}(s) \bigr)^+ ,
\\
\label{3.2.b} \Lambda(t) &=&  \max_{0 \le s \le t} \bigl(-r_2-g s
 +A(s)-\sigma V_{2}(s) \bigr)^+  .
\end{eqnarray}
These are modeled on (\ref{2.14}) and (\ref{2.18}), using the
identifications $  A (\cdot) \equiv L^Y (\cdot) $,\break $  \Lambda
(\cdot)
\equiv L^{R_2} (\cdot)$.

Such an approach is predicated on developing a theory for the unique
solvability of the system of equations (\ref{3.2.a}) and (\ref{3.2.b});
see Proposition \ref{Prop1} below and its proof in Section~\ref{sec5}.
The construction presented there expresses the continuous, increasing
and adapted processes $  A(\cdot) $, $  \Lambda(\cdot) $ as
\begin{equation}
\label{ProMeas} A(t)  =  \mathfrak{ A}   \bigl(t, (V_1 ,
V_2 ) \vert_{[0,t]}  \bigr) ,\qquad   \Lambda(t)  =
\mathfrak{ L}   \bigl(t, (V_1 ,  V_2 )
\vert_{[0,t]}   \bigr) ,\qquad   0 \le t < \infty ,
\end{equation}
progressively measurable functionals of the restriction $  (V_1 ,
V_2 )  \vert_{[0,t]}= \{  (V_1 (s),  V_2 (s)),   0 \le s \le t
\} $
of the planar Brownian motion $  (V_1 ,  V_2 ) $ on the
interval $  [0,t] $, and implies
\begin{equation}
\label{filt} \mathfrak{F}^{ (A, \Lambda)} (t)   \subseteq
\mathfrak{F}^{
(V_1, V_2)} (t) .
\end{equation}

\begin{proposition}
\label{Prop1}
Given the planar Brownian motion $  (V_1 (\cdot), V_{2}(\cdot)) $,
there exists a unique solution $  (A(\cdot), \Lambda(\cdot))  $ to
the system of equations (\ref{3.2.a}) and (\ref{3.2.b}); this is
expressible as a progressively measurable functional (\ref{ProMeas}).
\end{proposition}

\subsection{Constructing the gap, the laggard and the leader}
\label{sec4.1}

With the processes constructed so far, we introduce now the continuous
supermartingale
\begin{equation}
\label{3.3} Z(t)  :=  |y| - \lambda  t - \Lambda(t) + V^{ \flat} (t)
 , \qquad  0 \le t < \infty
\end{equation}
and its \textsc{Skorokhod} reflection at the origin
\begin{equation}
\label{3.4} G (t)  :=  Z(t) + 2  A(t) =  Z(t) + \max_{0 \le s \le t} \bigl(
-Z (s) \bigr)^+  \ge  0  , \qquad  0 \le t < \infty .
\end{equation}
(Here, the second equality is by virtue of (\ref{3.2.a}); the
nonnegative process $ G(\cdot) $ will play the r\^ole of the \textit{gap} between the leader and the laggard of the two semimartingales $
X_1 (\cdot) $, $ X_2 (\cdot) $ that we shall construct eventually,
in Section \ref{sec4.4}.)
We note that $   ( G (\cdot),   2  A(\cdot) ) $ is the
solution to
the \textsc{Skorokhod} reflection problem for the continuous
semimartingale $  Z(\cdot) $ of (\ref{3.3}), whose martingale part
has quadratic variation $  \langle V^\flat\rangle(t) = t $; and
with the help of the second property in (\ref{flat}) and of (\ref
{3.4}), (\ref{3.3}), we have the $ \mathbb{P}$-a.e. identities
\begin{equation}
\label{3.4.a} \int_0^\infty  \mathbf{
1}_{ \{ G (t) >0\}} \,  \mathrm{d} A (t)   = 0  ,\qquad   \int
_0^\infty  \mathbf{ 1}_{ \{ G (t) =0\}}
\,\mathrm{d} t  = 0 .
\end{equation}

\noindent$\bullet $
Let us introduce also the continuous semimartingale
\begin{equation}
\label{K} K(t)  :=  r_2 + g  t - A (t)+ \sigma
V_2(t)  ,\qquad   0 \le t < \infty
\end{equation}
and its \textsc{Skorokhod} reflection at the origin
\begin{equation}
\label{3.17.a} M (t)  :=   K (t) + \Lambda(t)   =   K(t) +
\max_{0 \le s \le t} \bigl( -K (s) \bigr)^+   \ge  0  .
\end{equation}
(Here, the second equality is by virtue of (\ref{3.2.b}), (\ref{K});
the nonnegative process $ M(\cdot) $ will play the r\^ole of the
\textit{laggard} of the two semimartingales $  X_1 (\cdot) $, $  X_2
(\cdot) $ that we shall construct in the next subsection.) The pair
$   ( M (\cdot),   \Lambda(\cdot) ) $ is the solution to
the \textsc{Skorokhod} reflection problem for the continuous
semimartingale $  K(\cdot) $ of (\ref{K}), whose martingale part
has quadratic variation $  \sigma^2 \langle V_2 \rangle(t) = \sigma^2 t $ with $  \sigma^2 >0 $, so we have the $  \mathbb{P}$-a.e.
identities
\begin{equation}
\label{3.18} \int_0^\infty  \mathbf{
1}_{\{M (t)>0\}}  \, \mathrm{d} \Lambda(t)  = 0  ,\qquad   \int
_0^\infty  \mathbf{ 1}_{ \{ M (t) =0\}}
\,\mathrm{d} t  = 0 .
\end{equation}
This last identity is a consequence of the second property in (\ref
{flat}), and of (\ref{3.17.a}), (\ref{K}).

\noindent$\bullet $
Finally, we introduce the continuous semimartingale
\begin{equation}
\label{3.16.a} N (t) :=  r_1 - h  t + A (t)+ \rho
V_1(t) ,\qquad   0 \le t < \infty
\end{equation}
by analogy with (\ref{2.16}), and note
\begin{equation}
\label{3.4.b} N(t) - M(t)  =  |y| - \lambda  t - \Lambda(t) +
V^{ \flat} (t) + 2  A(t)  =  G(t) \ge0  ,\qquad   0 \le t < \infty ,
\end{equation}
as well as the similarity of (\ref{3.16.a}) with (\ref{2.16}), and of
(\ref{3.17.a}) with (\ref{2.17}). The inequalities in (\ref
{3.17.a}), (\ref{3.4.b}) imply
\[
  \mathbb{P}  \bigl( N (t) \ge M (t) \ge0 ,   \forall  0 \le t < \infty
\bigr) =1 .
\]
Thus, the process $  N(\cdot) $ of (\ref{3.16.a}) is nonnegative;
it will play the r\^ole of the \textit{leader} of the two semimartingales
$  X_1 (\cdot) $, $  X_2 (\cdot) $ in the next subsection.

\noindent$\bullet $
Using results of Varadhan and Williams \cite{MR792398} and Reiman and Williams \cite{MR921820}
on Brownian motion with reflection in a wedge, we shall prove in Section
\ref{sec6} the following three propositions; related results have been
obtained by Burdzy and Marshall \cite{MR1231985,MR1227417}.

\begin{proposition}
\label{Prop3}
Under the condition $  1 / 2 \le\sigma^2 <1 $ of (\ref
{SigmaGeRho}), the planar process $   (N (\cdot) , M(\cdot))   $
with values in the acute (45-degree) wedge $   \mathfrak{ M} =\{ (n,
m) \in\mathbb{R}^2 \dvt 0 \le m \le n \}   $ never hits the corner of
the wedge:
\begin{equation}
\label{3.17.b} \mathbb{P}  \bigl( N (t)  >  0 ,  \forall   0 \le t <
\infty \bigr)  = 1 .
\end{equation}
\end{proposition}

\begin{proposition}
\label{Prop5}
In the case
\begin{equation}
\label{RhoGreSigma} 0   <   \sigma^2   <  1/2  ,
\end{equation}
the planar process $  (N (\cdot) , M(\cdot)) $ hits the corner of
the wedge $ \mathfrak{ M} $ with positive probability, that is, $
\mathbb{P}
  ( N (t) = 0,   \mbox{ for some }   t \in(0, \infty) ) >0$; this
probability is equal to one if, in addition,
$ g = h = 0 $.
\end{proposition}

\begin{proposition}
\label{Prop6}
With the processes $  G(\cdot) $ and $  M(\cdot) $ introduced in
(\ref{3.4}) and (\ref{3.17.a}), respectively, the solution $
(A(\cdot), \Lambda(\cdot))  $ of the system (\ref{3.2.a})--(\ref
{3.2.b}) satisfies the $\mathbb{P}$-a.e. identities of (\ref{3.4.a}),
(\ref
{3.18}) and
\begin{equation}
\label{ReimWill} \int_0^\infty\mathbf{
1}_{ \{ M(t) =0\} } \,  \mathrm{d}A(t) = 0  ,\qquad   \int_0^\infty
\mathbf{ 1}_{ \{ G(t) =0\} }   \,\mathrm{d} \Lambda(t) = 0 .
\end{equation}
\end{proposition}

\subsection{Unfolding the gap}
\label{sec4.1.a}

Theorem 1 in Prokaj \cite{MR2494646} guarantees that there exists an
enlargement $  ( \Omega, \mathfrak{F}, \mathbb{P}) $ of our filtered
probability space $  ( \widetilde{\Omega}, \widetilde{\mathfrak{F}},
\widetilde{\mathbb{P}}) $ with a measure-preserving mapping $  \pi
\dvtx  \Omega
\rightarrow\widetilde{\Omega} $; on this enlargement $  V_1(\cdot),
  V_2 (\cdot) $ are still independent Brownian motions, and there
exists a continuous semimartingale $  Y(\cdot) $ such that
\begin{equation}
\label{3.5} G (t)  =  \big| Y (t) \big| \quad  \mbox{and}\quad    Y(t)  =  y + \int
_0^t \overline{\operatorname{sgn}} \bigl( Y(s) \bigr)
\,  \mathrm{d} Z(s) ,\qquad   0 \le t < \infty .
\end{equation}
(The symmetric definition of the signum function, the first one in
(\ref{sgn}), is crucial here.) In other words, we represent the
\textsc{Skorokhod} reflection $  G(\cdot) $ of the semimartingale
$  Z(\cdot) $ in (\ref{3.4}), (\ref{3.3}) as the ``conventional''
reflection $  |Y(\cdot)| $ of an appropriate semimartingale $
Y(\cdot) $, related to $  Z(\cdot) $ via the \textsc{Tanaka}
equation in (\ref{3.5}). From this equation and (\ref{3.3}), we see
that the process $  Y(\cdot) $ satisfies the analogue of the
equation (\ref{2.11}):
\begin{equation}
\label{3.6} Y(t)  =  y + \int_0^t
\overline{\operatorname{sgn}} \bigl( Y(s) \bigr)   \bigl( - \lambda\, \mathrm{d} s -
\mathrm{d} \Lambda(s) + \mathrm {d} V^{ \flat} (s) \bigr) , \qquad  0 \le t <
\infty .
\end{equation}
Whereas, on the strength of (\ref{3.4.a}) and (\ref{ReimWill}), we
obtain also the $  \mathbb{P}$-a.e. properties
\begin{equation}
\label{3.4.aaa} \int_0^\infty  \mathbf{
1}_{ \{ Y (t) =0\}}  \, \mathrm{d} \Lambda (t)  = 0 ,\qquad   \int
_0^\infty  \mathbf{ 1}_{ \{ Y (t) =0\}}
\,\mathrm{d} t  = 0 .
\end{equation}

\noindent$\bullet $
In the interest of completeness, we review here this methodology from
Prokaj \cite{MR2494646} in the special case $  y=0 $: One considers the
zero set $  \mathfrak{Z} := \{ t \ge0 \dvt  G(t)=0 \} $ of the
continuous semimartingale $  G(\cdot) $ in (\ref{3.4}), and
enumerates as $  \{ \mathcal{E}_k \}_{k \in\mathbb{N}} $ the
components of the set $  [0,\infty) \setminus\mathfrak{Z} $, that is,
the excursions of $  G(\cdot) $ away from the origin. This is
carried out in a measurable manner, so that the event $  \{ t \in
\mathcal{E}_k \} $ belongs to the $  \sigma$-algebra $  \widetilde
{\mathfrak{F}} (\infty) := \sigma ( \bigcup_{0 \le\theta< \infty
} \widetilde{\mathfrak{F}} (\theta)  ) $ for all $  k \in
\mathbb{N} $ and $  t \ge0  $. Introducing independent Bernoulli
random variables $  \xi_1,   \xi_2, \ldots $ with $  \widetilde
{\mathbb{P}}  (\xi_k = \pm1 ) = 1/2 $, such that\vspace*{1pt} the sequence $
\{ \xi_k \}_{k \in\mathbb{N}} $ is independent of the $\sigma
$-algebra $  \widetilde{\mathfrak{F}}  (\infty) $, one sets
\[
\Phi(t)  :=  \sum_{k \in\mathbb{N}}   \xi_k
\mathbf{ 1}_{\{
  t  \in  \mathcal{E}_k \}}\quad   \mbox{and} \quad   \mathfrak{F} (t)  :=
\widetilde{\mathfrak{F}} (t) \vee\mathfrak {F}^{ \Phi} (t) ;\qquad  0 \le t
< \infty .
\]
The key observation from Prokaj \cite{MR2494646} is the balayage-type formula
\[
\Phi(\cdot)  G(\cdot)  =  \int_0^\cdot
\Phi(t) \, \mathrm{d} G(t)  =  \int_0^\cdot
\Phi(t) \, \mathrm{d} Z(t) ,
\]
with the second equality a consequence of (\ref{3.4.a}). Defining this
process above as $  Y (\cdot) := \Phi(\cdot)  G(\cdot)  $, one
observes the properties
$
| Y(\cdot)  |= G(\cdot) $, $    \overline{\mathrm
{sgn}}   (
Y(\cdot)  ) = \Phi(\cdot) $ and obtains the equation $
Y(\cdot) =
\int_0^\cdot  \overline{\operatorname{sgn}}   ( Y( t)  )\,\mathrm{d}
Z(t) $ from the equality of the first and third terms; this is (\ref{3.5})
for $  y=0 $.

\subsection{Constructing the various Brownian motions}
\label{sec4.2}

We are now in a position to trace the steps of the analysis we carried
out in Section \ref{sec3}, in reverse. Using the independent, standard
Brownian motions $  V_1 (\cdot),   V_2 (\cdot) $ we started this
section with, and the process $  Y(\cdot) $ we generated from them
in (\ref{3.5}), (\ref{3.6}) by enlarging the original probability
space, we introduce the two planar Brownian motions $   ( B_1
(\cdot),   B_2 (\cdot) ) $ and $   ( W_1 (\cdot),   W_2
(\cdot) ) $ via
\begin{eqnarray}
\label{3.7} B_1(t)  &:= & \int_0^t
\mathbf{ 1}_{ \{ Y (s) > 0 \} } \,  \mathrm {d} V_1 (s) + \int
_0^t \mathbf{ 1}_{ \{ Y (s) \le0 \} } \,  \mathrm{d}
V_2 (s)  ,
\\
\label{3.7.a} B_2(t)  &:=&   \int_0^t
\mathbf{ 1}_{ \{ Y (s) \le0 \} }  \, \mathrm{d} V_1 (s) + \int
_0^t \mathbf{ 1}_{ \{ Y (s) > 0 \} }  \, \mathrm{d}
V_2 (s)  ,
\end{eqnarray}
and
\begin{eqnarray}
 \label{3.8} W_1 (t)  &:=& \int_0^t
\Phi(s) \,  \mathrm{d} V_1 (s)  = \int_0^t
\overline{\operatorname{sgn}}  \bigl( Y (s) \bigr)  \, \mathrm{d} V_1 (s)
 = \int_0^t \operatorname{sgn} \bigl( Y (s) \bigr)
 \, \mathrm{d} V_1 (s)  ,
\\
\label{3.8.a} W_2 (t)  &:=& - \int_0^t
\Phi(s)  \, \mathrm{d} V_2 (s)  = -\int_0^t
\overline{\operatorname{sgn}}  \bigl( Y (s) \bigr) \,  \mathrm{d} V_2 (s)
 = -\int_0^t \operatorname{sgn}   \bigl( Y (s)
\bigr)  \, \mathrm{d} V_2 (s)  ,\quad
\end{eqnarray}
respectively. The last equalities in (\ref{3.8}), (\ref{3.8.a}) are
by virtue of the second equation in (\ref{3.4.aaa}), which also
implies that the equations of (\ref{2.3})--(\ref{2.6}) continue to hold.

Using these processes, we construct the Brownian motions $  W (\cdot)
$, $  V (\cdot)  $, $  Q (\cdot) $ and $  V^\flat(\cdot)   $
exactly as in (\ref{2.7}) and (\ref{2.9.b}); whereas by analogy with
(\ref{2.10}), and once again thanks to the second equation in (\ref
{3.4.aaa}), we note the representation
\begin{equation}
\label{BVW} W(\cdot) = \int_0^\cdot\overline{
\operatorname{sgn}}  \bigl( Y (t) \bigr)\, \mathrm{d} V^\flat(t) = \int
_0^\cdot\operatorname{sgn} \bigl( Y (t) \bigr)
\,\mathrm{d} V^\flat(t)  .
\end{equation}

\noindent$\bullet $ In conjunction with this last representation (\ref{BVW})
and the \textsc{Tanaka} formula (\ref{Tanaka2}), the properties of
(\ref{3.4.aaa}) allow us to write the equation (\ref{3.6}) for $
Y(\cdot) $ as
\begin{equation}
\label{3.6.b} Y (t)   =  y + \int_0^t
\operatorname{sgn} \bigl( Y(s) \bigr)  \bigl( - \lambda  \,\mathrm{d} s - \mathrm{d}
\Lambda(s) \bigr) + W(t)
\end{equation}
and obtain $ G(t)= | Y(t)| = |y| - \lambda  t - \Lambda(t) + V^{
\flat} (t) + 2  L^Y (t) = Z (t) + 2  L^Y (t)  ,   0 \le t < \infty
 $ on account of (\ref{3.3}). A comparison of this last expression
with (\ref{3.4}), using (\ref{2.19.c}) and (\ref{3.4.aaa}), gives
\begin{equation}
\label{3.17.d} 2  A(\cdot) =  2  L^Y (\cdot)  =
L^{|Y|} (\cdot) .
\end{equation}

\subsection{Naming the particles, then ranking them}
\label{sec4.4}

We can introduce now the continuous semimartingales
\begin{eqnarray}
\label{3.11} %
  X_1 (t) &:=&
 x_1 + \int_0^t ( g   \mathbf{
1}_{ \{ Y (s) \le
0 \}
} - h   \mathbf{ 1}_{ \{ Y (s) > 0 \} } )   \,\mathrm{d} s
\nonumber\\
&&{}+   \rho\int_0^t \mathbf{
1}_{ \{ Y (s) > 0 \} } \, \mathrm{d} B_1 (s) +  \int
_0^t \mathbf{ 1}_{ \{ Y (s) \le0 \} }  \bigl(
\sigma \, \mathrm{d} B_1 (s) + \mathrm{d} \Lambda(s) \bigr) ,
\\
&&  0
\le t < \infty,\nonumber   %
\\
\label{3.12} %
  X_2 (t) &:=&
x_2 + \int_0^t ( g   \mathbf{
1}_{ \{ Y (s) > 0 \} } - h   \mathbf{ 1}_{ \{ Y
(s) \le0
\} } )  \, \mathrm{d} s
\nonumber\\
&&{}+   \rho\int_0^t\mathbf{ 1}_{ \{ Y (s) \le0 \} }
  \,\mathrm{d} B_2 (s) +  \int_0^t
\mathbf{ 1}_{ \{ Y (s) > 0 \} }  \bigl( \sigma \, \mathrm{d} B_2 (s) +
\mathrm{d} \Lambda(s) \bigr) ,
\\
&& 0 \le t < \infty .\nonumber   %
\end{eqnarray}
From (\ref{2.5})--(\ref{2.6}), (\ref{3.1.a})--(\ref{3.1.b}), (\ref
{3.7})--(\ref{3.7.a}) and (\ref{3.4.aaa}), (\ref{3.6.b}), we obtain
for these two processes
\begin{eqnarray}
\label{3.13} X_1 (t) - X_2 (t)  &=&  y + \int
_0^t \operatorname{sgn}  \bigl( Y(s) \bigr)  \bigl(
- \lambda \, \mathrm{d} s - \mathrm{d}  \Lambda(s) + \mathrm{d}
V^{\flat}(s) \bigr) =  Y (t) ,
\\
\label{3.14} X_1 (t) + X_2 (t)  &=& \xi+ \nu t + V
(t) + \Lambda(t) , \qquad  0 \le t < \infty .
\end{eqnarray}

Repeating the analysis in Section \ref{sec3} and recalling the
notation of (\ref{2.22}), we obtain \textit{for all} $  t \in[0, \infty
) $ the analogues of (\ref{2.23b}), (\ref{2.24}), the skew representations
\begin{eqnarray}
\label{2.33} X_1 (t) &=& x_1 + \mu  t +
\rho^2 \bigl( Y^+ (t) - y^+ \bigr) - \sigma^2 \bigl( Y^-
(t) - y^- \bigr)
\nonumber\\[-8pt]\\[-8pt]
&&{} + \rho\sigma Q (t) + \bigl( \rho^2 -
\sigma^2 \bigr) \bigl( \tfrac{ 1 }{ 2 }  \Lambda(t) - A (t) \bigr) ,\nonumber
\\
\label{2.34}
 X_2 (t) &=& x_2 + \mu  t -
\sigma^2 \bigl( Y^+ (t) - y^+ \bigr) + \rho^2 \bigl( Y^-
(t) - y^- \bigr)
\nonumber\\[-8pt]\\[-8pt]
&&{} + \rho\sigma Q (t)+ \bigl( \rho^2 -
\sigma^2 \bigr) \bigl( \tfrac{ 1 }{ 2 }  \Lambda(t) - A (t) \bigr) .\nonumber
\end{eqnarray}

\noindent$\bullet $ Let us consider now as in (\ref{ranks}) the ranked versions
$  R_1 (\cdot) := X_1 (\cdot) \vee X_2 (\cdot)  $, $  R_2 (\cdot
) := X_1 (\cdot) \wedge X_2 (\cdot) $ of the semimartingales
introduced in (\ref{3.11}), (\ref{3.12}). From (\ref{3.13}), (\ref
{3.14}) we obtain
\begin{eqnarray}\label{3.15}
R_1 (t) + R_2 (t) &=&  X_1 (t)
+ X_2 (t)  =  \xi+ \nu t + V(t) + \Lambda(t) , \qquad  0 \le t <
\infty,
\nonumber\\[-8pt]\\[-8pt]
 R_1 (t) - R_2 (t) &=&   \big|
X_1 (t) - X_2 (t) \big| =   \big| Y (t) \big|  = G (t)  =  |y|
+ V^{ \flat} (t) - \lambda  t - \Lambda(t) + 2  A (t)\qquad\quad \nonumber
\end{eqnarray}
for their sum and gap, respectively, with the notation of Section
\ref{sec4.1}. These last two relations lead now with the help of (\ref
{3.16.a}), (\ref{3.17.a}) to the analogues of the equations (\ref
{2.16}) and (\ref{2.17}), that is,
\begin{eqnarray}
\label{3.16} R_1 (t) &\equiv&  N(t)  =  r_1 - h  t
+ \rho  V_1(t) + A (t) ,
\\
\label{3.17}  R_2 (t) &\equiv& M(t)  =  r_2 + g  t
+ \sigma  V_2(t)- A (t) + \Lambda(t)  .
\end{eqnarray}

\begin{rem}
\label{61}
Under the condition (\ref{SigmaGeRho}), the planar process $  (
X_1 (\cdot), X_2 (\cdot)  )$ constructed in (\ref{3.11}), (\ref
{3.12}) takes values in the punctured nonnegative quadrant $ \mathfrak
{ S} $. Indeed, Proposition \ref{Prop3} implies then
\begin{equation}
\label{3.17.f} \mathbb{P}  \bigl( X_1 (t) \vee X_2 (t)
 >  0 ,  \forall   0 \le t < \infty \bigr)  = 1 .
\end{equation}
Thus, the condition of (\ref{SigmaGeRho}) guarantees the absence of
``collisions at the origin''; see Ichiba and Karatzas \cite{MR2680554} and
Ichiba \textit{et~al.} \cite{2011arXiv1109.3823I} for similar conditions in a different context.
\end{rem}

\subsection{Denouement}
\label{sec4.6}

Let us start the final stretch of this synthesis by recalling the
second properties in each of (\ref{3.4.aaa}), (\ref{3.18}), which
lead to the $  \mathbb{P}$-a.e. identities
\begin{equation}
\label{3.19.b} \int_0^\infty\mathbf{
1}_{\{ X_1 (t) \wedge X_2 (t) =0\}} \, \mathrm {d} t  = 0 , \qquad  \int_0^\infty
\mathbf{ 1}_{\{ X_1 (t) = X_2
(t) \}} \, \mathrm{d} t  = 0 ;
\end{equation}
in particular, the ``nonstickiness'' conditions of (\ref{1.6}), both
along the boundary and along the diagonal, are satisfied. On the other
hand, the observation (\ref{2.19.a}) applied to the nonnegative
semimartingale $  R_2 (\cdot)  $ of (\ref{3.17}), together with the
properties of (\ref{3.19.b}), (\ref{3.18}) and (\ref{ReimWill}),
provides the characterization
\begin{eqnarray}
\label{3.20} %
  L^{X_1 \wedge X_2} (\cdot)  &=&
L^{R_2} (\cdot) =  \int_0^\cdot
\mathbf{ 1}_{\{ R_2 (t) =0\}} \, \mathrm{d} R_2 (t)
\nonumber\\[-8pt]\\[-8pt]
&=&  \int_0^\cdot\mathbf{ 1}_{\{ X_1 (t) \wedge X_2 (t) =0\}}
\bigl(  g \, \mathrm{d} t + \sigma \, \mathrm{d}V_2(t)- \mathrm{d} A
(t) + \mathrm{d}\Lambda(t) \bigr) =  \Lambda(\cdot) ; \nonumber  %
\end{eqnarray}
back into (\ref{3.4.aaa}), this gives the identity
\begin{equation}
\label{3.22} \int_0^\infty\mathbf{
1}_{\{ X_1 (t) = X_2 (t)\}}\,  \mathrm{d} L^{X_1 \wedge X_2} (t)  = 0 , \qquad  \mathbb{P}
\mbox{-a.e.}
\end{equation}

Arguing in a similar fashion, and applying the observation (\ref
{2.19.a}) to the nonnegative semimartingale $  G(\cdot) = R_1 (\cdot
) - R_2 (\cdot) \ge0 $ of (\ref{3.15}) in conjunction with the
properties (\ref{3.19.b})--(\ref{3.22}), (\ref{3.4.a}), (\ref
{3.4.aaa}) %
and (\ref{3.17.d}), we obtain by analogy with (\ref{2.19}) the $
\mathbb{P}
$-a.e. identity
\begin{eqnarray}
\label{3.21} %
 L^{R_1-R_2} (\cdot)  &=&  \int
_0^\cdot\mathbf{ 1}_{\{ X_1 (t) =
X_2 (t)\}} \, \mathrm{d}
\bigl( R_1 (t) - R_2 (t) \bigr)\nonumber
\\
 &=&  \int_0^\cdot\mathbf{ 1}_{\{ X_1 (t) = X_2 (t)\}}
\bigl( \mathrm{d} V^\flat(t) - \lambda \, \mathrm{d} t - \mathrm{d}
\Lambda(t) + 2 \,  \mathrm{d} A (t) \bigr)
\\
 &=& 2   A (\cdot)  = 2   L^Y (\cdot)  =  L^{|Y|} (
\cdot)  .\nonumber
\end{eqnarray}
With the help of (\ref{3.16}), (\ref{3.17}), we recover from these
last two observations the equations (\ref{2.20}), (\ref{2.21}) for
the ranks; whereas from (\ref{ReimWill}) and the identifications $
A(\cdot) = L^Y (\cdot) = L^{|Y|} (\cdot) / 2 $ in (\ref{3.21}), we
deduce the $ \mathbb{P}$-a.e. identities
\begin{equation}
\label{3.23} \int_0^\infty\mathbf{
1}_{\{ X_1 (t) \wedge X_2 (t)= 0\}}\,  \mathrm {d} L^{X_1 - X_2} (t) = 0  ,\qquad   \int
_0^\infty\mathbf{ 1}_{\{
X_1 (t) \wedge X_2 (t)= 0\}} \, \mathrm{d}
L^{X_2 - X_1} (t) = 0 .
\end{equation}
We conclude from (\ref{3.17.f})--(\ref{3.23}) and Proposition \ref
{Prop3} that we have proved the following result.

\begin{proposition}
\label{Prop2}
The continuous, planar semimartingale $  \mathcal{X} (\cdot) = (
X_1 (\cdot),   X_2 (\cdot)  )'  $ defined in (\ref{3.11}),
(\ref{3.12}) takes values in the quadrant $  [0, \infty)^2 $. It
satisfies the dynamics of (\ref{1.3}) and (\ref{1.4}); the properties
(\ref{1.6}), (\ref{1.7}) and (\ref{Ouk3}); as well as the
representations (\ref{2.23b}), (\ref{2.24}).

Under the condition (\ref{SigmaGeRho}), the planar semimartingale $
\mathcal{X} (\cdot)  $ takes values in the punctured nonnegative
quadrant of (\ref{1.2}), that is, never hits the corner of the quadrant.
\end{proposition}

\begin{rem}
\label{31}
From the equations (\ref{3.16}), (\ref{3.17}) and (\ref{ProMeas}),
we express the continuous, adapted processes $  R_1(\cdot) $, $
R_2 (\cdot) $ as progressively measurable functionals
\[
R_1 (t)  =  \mathfrak{R}_1   \bigl(t,
(V_1 ,  V_2 ) \vert_{[0,t]}  \bigr) ,\qquad
R_2 (t)  =  \mathfrak{R}_2   \bigl(t,
(V_1 ,  V_2 ) \vert_{[0,t]}   \bigr) ,\qquad  0
\le t < \infty
\]
of the restriction $  (V_1 ,  V_2 )  \vert_{[0,t]}= \{  (V_1
(s),  V_2 (s)),  0 \le s \le t  \} $ of the planar Brownian motion
$  (V_1 ,  V_2 ) $ on the interval $  [0,t] $; this implies
\begin{equation}
\label{flit2} \mathfrak{F}^{ (R_1, R_2)} (t)   \subseteq
\mathfrak{F}^{ (V_1,
V_2)} (t)  .
\end{equation}
From these observations, from the identifications $  L^{R_2} (\cdot)
\equiv\Lambda(\cdot) $ and $  L^Y (\cdot) \equiv A(\cdot) $ in
(\ref{3.20}), (\ref{3.21}), and from the analysis of Section \ref
{sec3} that culminates with the equations (\ref{2.16})--(\ref{2.17}),
we deduce that \textit{the distribution of the vector of ranked processes
$  ( R_1 (\cdot), R_2 (\cdot) ) \equiv ( N (\cdot), M
(\cdot) )  $ in \textup{(\ref{ranks})} is determined uniquely.}

A more elaborate analysis, carried out in Section \ref{sec7}, will
show that
uniqueness in distribution (indeed, pathwise uniqueness up until its first
visit to the corner of the quadrant) holds also for the vector process
$
\mathcal{X}(\cdot)=  ( X_1 (\cdot), X_2 (\cdot) )'  $ of the
particles' positions by name in (\ref{3.11}), (\ref{3.12}). If the two
particles never collide with each other at the origin, as is the case under
condition (\ref{SigmaGeRho}), then pathwise uniqueness -- thus also
uniqueness in distribution -- holds for all times.
\end{rem}

\begin{rem}
\label{51}
Coupled with (\ref{flit}) of Remark \ref{22}, the filtration
inclusion in (\ref{flit2}) gives in the nondegenerate case $  \sigma
<1 $ the identity
\[
\mathfrak{F}^{ (R_1, R_2)} (t)   =   \mathfrak{F}^{ (V_1, V_2)} (t) ,\qquad
 0 \le t < \infty .
\]
\end{rem}

\begin{figure}

\includegraphics{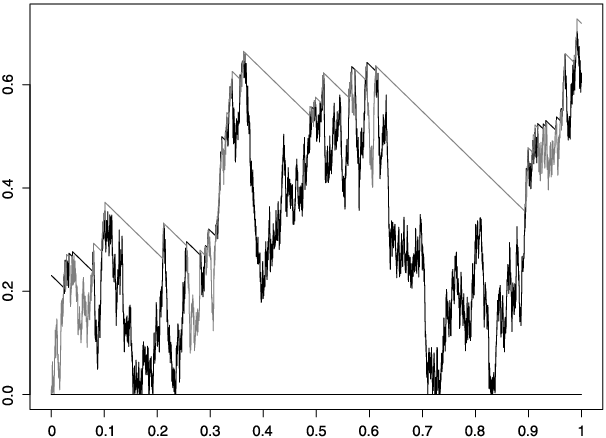}

\caption{Simulated processes; Black $ =X_1(\cdot) $, Gray $ =X_2
(\cdot) $.}\vspace*{-2pt}\label{fig1}
\end{figure}

The Figure \ref{fig1} is a simulation of the processes $ X_1(\cdot)
 $ and $ X_2(\cdot)  $ in the degenerate case
with \mbox{$  \rho=0 $}, $  \sigma=1 $ and $  g = h =1 $, taken from
Fernholz \cite{Fernholz2011a}; we are grateful to \textsc{Dr. Robert Fernholz}
for granting us permission to reproduce it here. The figure depicts
clearly also the laggard (from (\ref{3.17})) and the leader (from
(\ref{3.16})) of the two processes -- the latter as a function of
finite first variation which ascends by the continuous but singular,
local time component $  L^Y (\cdot) = A(\cdot) $ and descends by
straight line segments (``ballistic motion'') of slope $-1$.


\section{\texorpdfstring{Proof of Proposition \protect\ref{Prop1}}{Proof of Proposition 4.1}}
\label{sec5}

Given the planar Brownian motion $  (V_1 (\cdot), V_{2}(\cdot))  $,
we shall
apply the idea of the proof of Theorem~1 of Harrison and Reiman \cite{MR606992} to study a system
of equations for
$(2   A(\cdot), \sqrt{2}   \Lambda(\cdot))^{\prime}   $
equivalent to
(\ref{3.2.a})--(\ref{3.2.b}):
\begin{eqnarray}
\label{eq: 6 A-a} %
2   A(t) & =& \max_{0\le s \le t}
\bigl({-   \vert y \vert+   \lambda  s +   \Lambda(s) -
V^{\flat}(s) } \bigr)^+  ,
\nonumber\\
\sqrt{2}  \Lambda(t) &=& \max_{0\le s \le t} \bigl({-   \sqrt{2}
r_{2} - \sqrt{2}  g  s + \sqrt{2}  A (s) - \sqrt{2}   \sigma
V_2(s) } \bigr)^+  ,
\\
 &&0 \le t < \infty ,\nonumber
\end{eqnarray}
with $ V^{ \mathbf{ \flat}}( \cdot) = \rho  V_1 (\cdot) - \sigma
  V_2 (\cdot)  $ as in (\ref{2.9.b}). Schematically, we shall write
$ {\mathbf y} (\cdot) = \bolds{\pi}_{\mathfrak w}
({\mathbf y})(\cdot) $ for this system, with the two-dimensional
\textsc{Skorokhod} map $ \bolds{\pi}_{\mathfrak w} $ introduced in
Harrison and Reiman \cite{MR606992}:
\begin{equation}
\label{Pi} C_0\bigl([0, \infty); \mathbb{R}^{2}\bigr)
\ni{\mathbf y} (\cdot)   \mapsto   \bolds{\pi}_{\mathfrak w} ({\mathbf
y}) (\cdot) := \sup_{0 \le
s \le  \cdot} \bigl[  {\mathbf H}  {\mathbf y}(s) - {
\mathfrak w}(s)  \bigr]^{+} \in C_0\bigl([0, \infty
); \mathbb{R}^{2}\bigr)  .
\end{equation}
Here the subscript indicates the pinning $  {\mathrm y}_1 (0) =
{\mathrm y}_2 (0) = 0  $; the matrix $  {\mathbf H}  $ is defined as
\[
{\mathbf H} = \frac{1}{\sqrt{2}} %
\pmatrix{ 0 & 1
\cr
1 & 0 } %
 ,
\]
and the component processes of the vector $  {\mathfrak w}(\cdot) =
({\mathfrak w}_{1} (\cdot), {\mathfrak w}_{2}(\cdot))^{\prime}  $ are
\begin{eqnarray}
\label{eq: w-a} %
  {\mathfrak w}_1(t) & =&
\llvert {y}\rrvert - \lambda  t + V^{\flat}(t) = \llvert {y}\rrvert -
\lambda  t + \rho  V_1(t) - \sigma V_2 (t)   ,\nonumber
\\
{\mathfrak w}_2 (t) &=& \sqrt{2}   r_2 + \sqrt{2}  g
t + \sqrt {2}  \sigma  V_2(t)  ,
\\
&& 0
\le t < \infty .\nonumber
\end{eqnarray}
The supremum and the positive part $  ( \cdot )^{+}:= \max(
\cdot , 0) $ are taken for each element in the vectors.

The matrix $  {\mathbf H}  $ has spectral radius $  1 / \sqrt{2}
$, so the mapping $  \bolds{\pi}_{\mathfrak w}  $ is a continuous,
contraction mapping (Theorem 1 of Harrison and Reiman \cite{MR606992}); in particular, in
terms of
the sup-norm
\[
  \Vert{\mathbf y} \Vert_{T}  :=  \max_{1 \le i \le2}
\sup_{0
\le s \le T}   \big\vert  \mathrm{y}_i (s) \big\vert
\]
for every $  {\mathbf y},  {\mathbf y}^{\flat} \in C_0 ([0, \infty
); \mathbb{R}^2) $ and $  T \in[ 0, \infty)   $ we have the bounds
\begin{eqnarray*}
  && \bigl\llvert {\bolds{\pi}_{\mathfrak w} ({
\mathbf y}) - \bolds{\pi}_{\mathfrak w} \bigl({\mathbf y}^\flat
\bigr)}\bigr\rrvert_{T}
\\
&&\quad  = \max\lleft[ \sup_{0 \le s \le T} \bigg\vert
\sup_{0 \le u \le
s} \biggl( - {\mathfrak w}_1 (s) +
\frac{1}{\sqrt{2}} \mathrm{y}_2(s) \biggr)^+ - \sup_{0 \le u \le s}
\biggl( - {\mathfrak w}_1 (s) + \frac
{1}{\sqrt{2}}
\mathrm{y}_2^\flat (s) \biggr)^+ \bigg\vert  ,\rright.
\\
&&\hspace*{46pt}\lleft.\sup_{0 \le s \le T} \bigg\vert\sup_{0 \le u \le s} \biggl( - {\mathfrak
w}_2 (s) + \frac{1}{\sqrt{2}} \mathrm{y}_1(s) \biggr)^+ -
\sup_{0 \le u \le s} \biggl( - {\mathfrak w}_2 (s) +
\frac{1}{\sqrt {2}} \mathrm{y}_1^\flat(s) \biggr)^+ \bigg\vert
 \rright]
\\
&&\quad  \le\max\lleft[ %
 \sup_{0 \le s \le T}
\bigg\vert{-} \mathfrak {w}_1 (s) + \frac
{1}{\sqrt{2}} \mathrm{y}_2(s)
- \biggl( - {\mathfrak w}_1(s) + \frac{1}{\sqrt{2}}
\mathrm{y}_2^\flat(s) \biggr) \bigg\vert,\rright.
\\
&&\hspace*{46pt}\lleft.
 \sup_{0 \le s \le T} \bigg\vert{-} \mathfrak {w}_2 (s) +
\frac{1}{\sqrt{2}} \mathrm{y}_1(s) - \biggl( - {\mathfrak
w}_2(s) + \frac{1}{\sqrt{2}} \mathrm{y}_1^\flat(s)
\biggr) \bigg\vert
\rright]
\\
&&\quad   \le  \frac{1}{ \sqrt{2} }  \big\Vert{
\mathbf y} - {\mathbf y}^\flat \big\Vert_T  .
\end{eqnarray*}
We have used here the contraction property %
of the maximum: for two continuous, real-valued functions $  \mathrm
{y} (\cdot),
 \mathrm{y}^\flat(\cdot)  $ we have
\begin{eqnarray*}
\Bigl\llvert {\sup_{0 \le s \le T} \mathrm{y} (s) - \sup_{0 \le s \le T}
\mathrm{y}^\flat(s) }\Bigr\rrvert &\le&\sup_{0 \le s \le T}  \bigl\llvert
{\mathrm{y} (s) - \mathrm {y}^\flat(s)}\bigr\rrvert  , \qquad  0 \le T <
\infty  ,
\\
\bigl\llvert {\bigl( \mathrm{y}(s)\bigr)^+ - \bigl(\mathrm{y}^\flat(s)
\bigr)^+ }\bigr\rrvert &\le& \bigl\llvert {\mathrm{y}(s) - \mathrm{y}^\flat(s)}
\bigr\rrvert  ,\qquad  0 \le s < \infty  .
\end{eqnarray*}
It follows from this contraction property that, given $  {\mathfrak
w}(\cdot)  $, the solution $  {\mathbf y}(\cdot)   $ to the system of
(\ref{eq: 6 A-a}) can be obtained as the \textit{unique} limit of a standard
\textsc{Picard--Lindel\"of} iteration: starting with $  {\mathbf
y}^{(1)}(\cdot) \equiv0 $, iterating $  {\mathbf y}^{(n+1)}:=
\bolds{\pi}_{\mathfrak w} (  {\mathbf y}^{(n)})  $ for $  n = 1,
2,\dots ,$ we obtain the uniform convergence on compact intervals $
\lim_{n \to\infty} \Vert {\mathbf y} - {\mathbf y}^{(n)}\Vert_{T} =
0 $ for $  T \in( 0, \infty)  $.
This concludes the proof of Proposition \ref{Prop1}.

It might be interesting to investigate possible connections between this
construction and the \textsc{Skorokhod} map on an interval, studied by
Kruk \textit{et al.} \cite{MR2349573}.

\section{\texorpdfstring{Proof of Propositions \protect\ref{Prop3}--\protect\ref{Prop6}}
{Proof of Propositions 4.2--4.4}}
\label{sec6}

Let us recall from (\ref{3.16.a}), (\ref{3.17.a}) that the process $
(N(\cdot),M(\cdot)) $ %
is a reflected planar Brownian motion in the $45$-degree wedge $
\mathfrak M $ with orthogonal reflection on the faces:
\begin{eqnarray}
\label{R.1} N (t) &=&  r_1 - h  t + \rho V_1(t) + A
(t) , \qquad 0 \le t < \infty ,
\\
\label{R.2} M (t) &=&  r_2 + g  t + \sigma V_2 (t)- A
(t) + \Lambda(t) , \qquad  0 \le t < \infty .
\end{eqnarray}
This process does not hit the corner of the wedge $ \mathfrak{M} $,
if and only if it does not reach the corner during the time-horizon $
[0, T] $ for any $ T \in(0, \infty) $. Hence, in the
nondegenerate case we can assume that the drift coefficients $ g, h
$ are equal to zero. Indeed, after a suitable change of probability
measure on $ \mathfrak F(T) $, under the new measure and on the
finite time-horizon $  [0, T] $ the process $   (
V_1(t)-(h/\rho)  t ,V_2(t)+(g/\sigma) t  )_{t\in[0,T]} $
becomes then a two-dimensional standard Brownian motion, that is, a
planar Brownian motion without drift.

Thus, in what follows we shall assume $ h=g=0 $ and apply the transformation
$ \mathfrak{T}=\operatorname{diag}(1/\rho , 1/\sigma) $ to the
process $
(N(\cdot),M(\cdot)) $ as in (\ref{R.1}), (\ref{R.2}). In matrix
form, we can write
\begin{equation}
\label{eq: plRBM} \mathfrak{T} %
\left[\matrix{ N(t)
\cr
M(t) }\right] %
= \mathfrak{T} %
\left[\matrix{
r_1
\cr
r_2 }\right] %
+ %
\left[\matrix{
V_1(t)
\cr
V_2(t)
} \right]%
+ \mathfrak{T}   %
\left[\matrix{ 1 & 0
\cr
-1 & 1 }\right] %
\left[\matrix{ A(t)
\cr
\Lambda(t)
}\right] %
  ;\qquad   0 \le t < \infty .
\end{equation}
Then the transformed process is a reflected Brownian motion in the wedge
$\mathfrak{W}=\mathfrak{T}(\mathfrak{M})$ %
with oblique reflection on the faces.

We need to compute the angle $\bolds{\xi}$ of the wedge $ \mathfrak{W}
$ at the corner, and also the reflection angles
$  \bolds{\theta}_i \in(- \uppi  /   2, \uppi  /   2 )   $
$  i = 1, 2  $ measured from the inward normal vector on the
boundary, and
positive if and only if they direct the process toward the corner; see
Figure \ref{fig:2} below. We introduce also the scalar parameter
\begin{equation}
\label{eq: alpha} \bolds{\alpha}   :=   (\bolds{\theta}_1 + \bolds{\theta}_2)   /   \bolds{\xi} .
\end{equation}

According to Theorems 2.2 and 3.10 of
Varadhan and Williams \cite{MR792398}, the reflected Brownian motion in the wedge $\mathfrak{W}$
never hits the corner of the wedge $  \mathfrak W   $, with
probability one,
if $  \bolds{\alpha} \le0  $; hits the corner with probability one,
if $
\bolds{\alpha} >0   $; and is well-defined by the corresponding submartingale
problem for all times, starting at any initial point \textit{including the
corner}, if $  \bolds{\alpha} < 2  $.

Now the faces of the wedge $ \mathfrak{W} $ are given by half-lines
emanating from the origin and parallel to the vectors $ \mathfrak
{v}_1 =
\mathfrak{T}(1,0)'=(1/\rho, 0)' $ and $ \mathfrak{v}_2 =\mathfrak{T}
(1,1)'=(1/\rho ,1/\sigma)' $, hence %
\begin{equation}
\label{eq: XI} \cos(\bolds{\xi}) =  \frac{ \langle{\mathfrak{v}_1},{\mathfrak{v}_2} \rangle} {
\llVert {\mathfrak{v}_1}\rrVert \llVert {\mathfrak{v}_2}\rrVert }  =  \frac{\rho^{-2}} { \sqrt{\rho^{-2}+\sigma^{-2} }\sqrt{\rho^{-2} } } =
  \sigma .
\end{equation}
The reflection vector on the face of the wedge $ \mathfrak{W} $
parallel to $  \mathfrak{v}_2 =(1/\rho ,1/\sigma)' $ is $
\bolds{\nu}_2 :=   (1/\rho , -1/\sigma  )' $, while the normal
vector of this face pointing
inward is $ \mathbf{n}_2 =(\rho,- \sigma)' $. Then
\[
\bolds{\nu}_2-\frac{ \langle{\bolds{\nu}_2},{\mathbf{n}_2}
\rangle}{\llVert {\mathbf{n}_2}\rrVert^2} \mathbf{n}_2  =  \bolds{
\nu}_2-2  \mathbf{n}_2 =   (1/\rho-2\rho ,  2\sigma-1/
\sigma )' =  %
\bigl(\sigma^2-
\rho^2\bigr)\cdot(1/\rho , 1/\sigma )' ,
\]
so $ \bolds{\nu}_2 $ points towards to the corner exactly when $
\sigma^2<\rho^2 $. That is, $  \bolds{\theta}_2>0  $ holds if and only if $
\sigma^2<1/2 $.

The other face of the wedge $ \mathfrak{W} $ is parallel to $
(1,0)' $ and the reflection vector on this face is
$(0, 1/\sigma)'$, which is a normal vector to this face, whence $ \bolds{\theta}_1=0 $.

\begin{pf*}{Proof of Proposition \ref{Prop3}}
 If $  1> \sigma^{2} \ge1
/   2  $, then %
\begin{equation}
\label{eq: neg-alpha} \bolds{\theta}_1 =0  ,\qquad   \bolds{\theta}_2 \leq0  ,\qquad  \bolds{\xi } > 0 ,
\end{equation}
so $  \bolds{\alpha} \le0 $. Thus, according to the result of Harrison and Reiman \cite{MR606992} and %
Varadhan and Williams \cite{MR792398}, with probability one the process $ (N(\cdot
),M(\cdot)) $ never hits the corner of the wedge $ \mathfrak{M} $.
The case $  \sigma^2 =1 $ is discussed in Section \ref{6.111}.
\end{pf*}

\begin{pf*}{Proof of Proposition \ref{Prop5}}
In the case $  0 < \sigma^{2} < 1  /   2  $, we obtain $  \cos
(\bolds{\xi}) < 1/\sqrt{2}$ from (\ref{eq: XI}) and
\[
\bolds{\theta}_1 = 0  ,\qquad  0< \bolds{\theta}_2 < \uppi/ 2  ,\qquad
  \bolds{\xi} > \uppi  /   4  ,
\]
so $  0 < \bolds{\alpha}<2 $. Then the process $ (N(\cdot),M(\cdot
))$ hits the corner of the wedge $ \mathfrak{M} $ almost surely.
This gives the result for
$g=h=0$.

When $ g+h>0 $, we can only ascertain that the process $ (N(\cdot
),M(\cdot)) $ hits the corner of the wedge with positive probability,
due to the measure change step that we deployed to reduce the general
case to the driftless case.
\end{pf*}

\subsection{\texorpdfstring{Nonattainability of the corner in the degenerate case $\rho=0$}
{Nonattainability of the corner in the degenerate case rho = 0}}
\label{6.111}

In this subsection, we develop the proof of Proposition \ref{Prop3} in
the degenerate case $ \sigma=1, \rho=0$. The equations of \eqref
{R.1}, \eqref{R.2} for the ranked processes $ N (\cdot) \ge M (\cdot
) \ge0  $ simplify then to
\begin{eqnarray}
\label{R1} N (t) &= & r_1 - h  t + A (t) , \qquad 0 \le t < \infty,
\\
\label{R2} M (t) &=&  r_2 + g  t +   V (t)- A (t) + \Lambda(t) ,
\qquad  0 \le t < \infty ;
\end{eqnarray}
we recall that the ``regulating'' continuous, increasing and adapted
processes $  A(\cdot),  \Lambda(\cdot)  $ of (\ref{3.2.a}),
(\ref{3.2.b}) satisfy the $  \mathbb{P}$-a.e. requirements
\begin{equation}
\label{R3} \int_0^\infty  \mathbf{
1}_{ \{ N (t) > M(t) \}}  \, \mathrm{d} A (t)  = 0  ,\qquad   \int_0^\infty
  \mathbf{ 1}_{ \{ M (t) > 0\}}  \, \mathrm{d} \Lambda(t)  = 0
\end{equation}
as in (\ref{3.4.a}), (\ref{3.18}). We introduce the stopping time $
\tau=\inf \{{t>0 \dvt N(t)=0} \} $, and shall establish
below the
property $\mathbb{P}(\tau=\infty)=1$.

For this, it suffices to consider $  h>0 $; for if $  h=0 $, the
process $ N (\cdot) $ is nondecreasing, and there is nothing to
prove. Whereas, by the Girsanov theorem, it is enough to deal with the
case $g=0$. We shall present two distinct, very different arguments.

\noindent$\bullet $ \textsc{First argument:} In the manner of Section \ref
{sec4.1.a}, we consider the unfolded process $ (N (\cdot),$
$M^\dagger(\cdot)) $, where $ M (\cdot)=\llvert {M^\dagger(\cdot
)}\rrvert
$ and $ \mathrm{d}M^\dagger(t)=\sgn*(M^\dagger(t)) \, \mathrm{d}M
(t) $, stopped upon reaching the corner:
\begin{eqnarray*}
N (t)= r_1 - h  t +A (t) ,\qquad  M^\dagger(t)=
r_2 + V^\dagger(t)- \int_0^t
\sgn*\bigl(M^\dagger(s)\bigr) \, \mathrm{d}A(s)
\end{eqnarray*}
for $ 0 \le t\leq\tau $; and $ M (t)=N (t)=0 $ for $ t\geq\tau
 $.
Here $  V^\dagger(\cdot) :=\int_0^\cdot\sgn*(M^\dagger(s))
\,\mathrm{d}V(s) = \int_0^\cdot\sgn(M^\dagger(s))\,  \mathrm
{d}V(s)  $ is of course standard Brownian motion; the equality of the
two stochastic integrals follows by applying the second property of
(\ref{flat}) to the semimartingale $  M(\cdot) $ of (\ref{R2}).

The planar process $ (N(\cdot),M^\dagger(\cdot)) $ evolves in the
cone $  \{{(x,y) \dvt 0\leq\llvert {y}\rrvert \leq x} \} $,
with normal reflection
on the faces and absorption when the corner of the cone is reached. It
is clear that $(N (\cdot),M^\dagger(\cdot))$ is a strong Markov
process; when started at $(x,y)$, the distribution of this process will
be denoted by $\mathbb{P}_{(x,y)}$. We shall show that $ \mathbb{P}(\tau=\infty)=1 $.

For this purpose, we define a Markov chain $   \{ ( \widetilde
{S}_n,Y_n)  \}_{n\geq0} $ with
\[
\widetilde{S}_n := \log_{(2)} \bigl( N(
\tau_n) / r_1 \bigr)  ,\qquad    Y_n :=
 M^\dagger(\tau_n)/N(\tau_n) ,
\]
where we have set $ \tau_0:=0 $ and recursively
$ \tau_{n+1}:=\inf \{{t>\tau_n \dvt (N (t)/N(\tau_n))\notin
[1/2,2]} \} $, and noted that $ \tau_n $ is a.s. finite,
for all $
n\in\mathbb{N}_0 $. The state-space of this Markov chain is $\mathbb{Z}
\times[-1,1]$.
For $ k>0 $, we shall denote
\[
\widetilde{\rho}_k := \inf  \{   n \in\mathbb{N}_0
 \dvt  \widetilde{S}_n = \widetilde{S}_0 + k \} .
\]
A simple sufficient condition for the nonattainability of the origin
by $ N(\cdot) $ is that
\begin{equation}
\label{eq:SC} \mbox{$\widetilde{\rho}_k $ is finite almost surely
for all $ k>0$.}
\end{equation}
Indeed, on $ \{{\tau<\infty} \}$ the sample path of the process
$ N(\cdot) $ is bounded, because it is continuous. It is clear from
this that
$  \{{\tau< \infty} \}\subset\bigcup_{k} \{
{\widetilde{\rho }_k=\infty} \} $,
and $ \mathbb{P}(\tau<\infty)=0 $ follows from~\eqref{eq:SC}.

We compare $  \widetilde{\mathbf{ S}}=  \{ \widetilde{S}_n
\}_{n\geq0} $ with a random walk $  \mathbf{ S}=  \{ S_n  \}_{n\geq0} $ on $ \mathbb{Z} $, which starts at $S_0=0$ and is
defined by
\[
\mathbb{P}({S_{n+1}=S_n+1 \mid S_0,
\dots,S_n})= \mathfrak{p}\bigl(r_1 2^{S_n}
\bigr),\qquad   \mathbb{P}({S_{n+1}=S_n-1 \mid S_0,
\dots,S_n})=1- \mathfrak{p}\bigl(r_1 2^{S_n}
\bigr)
\]
as its transition probabilities, where
\[
\mathfrak{p} (r) :=\inf_{y\in[-r,r]}\mathbb{P}_{(r,y)}\bigl({N(
\tau_1)=2r}\bigr).
\]
To wit: we start the process $(N(\cdot),M^\dagger(\cdot))$ from a
point on the vertical
line at $ x=r $, and denote by $  \mathfrak{p}(r) $ the greatest
lower bound on the probability that $ N(\cdot) $ hits $ 2r $
before it hits $ r/2 $. We prove below that
\begin{equation}
\label{eq:p} \lim_{r\to0+} \mathfrak{p}(r)=1\quad  \mbox{and}\quad
\mathfrak {p}(r)>0\qquad  \mbox{for all $r>0$}.
\end{equation}
On a suitable extension of our probability space there is a coupling between
$  \widetilde{\mathbf{ S}} $ and the random walk $  \mathbf{ S}
$, such that $ \widetilde{S}_n\geq S_n $, $ \forall   n \in
\mathbb{N}_0 $. Therefore, the sufficient condition \eqref{eq:SC}
will be established as soon as we show that $ \rho_k=\inf \{
{n\geq 0 \dvt S_n=k} \}\geq\tilde\rho_k $ is a.s. finite, for
all $ k>0 $.

Let $  \ell_0\leq0 $ be such that $  \mathfrak{p}(r_12^\ell)\geq
1/2 $ holds for $  \ell\leq\ell_0 $, and consider a simple,
symmetric random walk $  \widehat{\mathbf{S}}= \{\widehat
{S}_n \}_{n\geq0} $ on the state space $  \{{\ell \dvt \ell
\leq \ell_0} \} $, with $\ell_0$ as both reflecting barrier
and starting
point. Using coupling again, we can assume that $  \widehat{S}_n\leq
S_n $ holds for all $ n \in\mathbb{N}_0 $. Fix $ k>0 $ and
recall that $ \widehat{\mathbf{S}} $ is recurrent. Thus, if $
\mathbf{S} $ does not reach $ k  $ before reaching $  \ell_0-1
$, it will return to the level $  \ell_0 $ almost surely; and from
$  \ell_0 $ it will reach $ k $ before reaching $  \ell_0-1 $,
with some positive probability. If, on the other hand, $ \mathbf{S}
$ hits $  \ell_0-1 $ first, then the whole thing starts again and
finally $ k $ is reached almost surely by $\mathbf{S}$. By a
standard renewal argument, this implies that $ \mathbf{ S} $ reaches
$S_0+k$ almost surely, and that
$\rho_k=\inf \{{n\geq0 \dvt S_n=S_0+k} \}$ is finite almost surely.

It remains now only to argue \eqref{eq:p}. Let $r>0$ and observe that,
if $ N(0)=r $ and $ \llvert {M^\dagger(0)}\rrvert \leq r $, then $
N(\tau_1)=r/2 $ implies that
$  \tau_1\geq(r/2h)$ and $ \llvert {V^\dagger(r/2h)}\rrvert <4r $, hence
\[
1- \mathfrak{p}(r) =  \sup_{y\in[-r,r]}\mathbb{P}_{(r,y)}\bigl(N(
\tau_1) =  r/2\bigr)  \leq  \mathbb{P}\bigg(\biggl\llvert {V^\dagger
\biggl({\frac{r}{2h}} \biggr)}\biggr\rrvert <4r\bigg),
\]
that is,
\[
\mathfrak{p}(r)\geq1-\mathbb{P}\bigg({\biggl\llvert {V^\dagger \biggl({
\frac
{r}{2h}} \biggr)}\biggr\rrvert <4r}\bigg) = \mathbb{P}\bigg( { \bigl\llvert
{V^\dagger(1)}\bigr\rrvert >4\sqrt{2h r } }\bigg).
\]
This justifies \eqref{eq:p} and completes the proof of the
nonattainability of
the corner in the degenerate case, thus also the proof of Proposition
\ref{Prop3}.

\noindent$\bullet $ \textsc{Second argument}: Here follows another argument,
due to \textsc{Dr.$ $E.$ $Robert Fernholz}; we shall take $  r_2=0
$ for simplicity. With $  B(\cdot) $ a standard Brownian motion, we
denote by $  \Gamma(\cdot) $ the \textsc{Skorokhod} reflection of
the process $  r_1- h t - B(t) ,  0 \le t < \infty $, and by $
\Delta(\cdot) $ the \textsc{Skorokhod} reflection of $  r_1- h t
+ B(t) ,  0 \le t < \infty $:
\begin{equation}
\label{6.1} \Gamma(t)  =  r_1- h t - B(t) + L^{\Gamma} (t)
  \ge  0 , \qquad  \Delta(t)  =  r_1- h t + B(t) +
L^{\Delta} (t)   \ge  0 ,
\end{equation}
where the continuous, increasing processes
\[
L^{\Gamma} (t) :=  \max_{0 \le s \le t} \bigl( h  s + B(s)
-r_1 \bigr)^+ ,\qquad    L^{\Delta} (t) :=
\max_{0 \le s \le t} \bigl( h  s - B(s) -r_1 \bigr)^+
\]
satisfy the $  \mathbb{P}$-a.s. identities $  \int_0^\infty\mathbf
{ 1}_{\{
\Gamma(t) >0\}}  \, \mathrm{d} L^{\Gamma} (t)=0 $, $   \int_0^\infty\mathbf{ 1}_{\{ \Delta(t) >0 \} } \,  \mathrm{d} L^{
\Delta} (t)=0 $. We define
\begin{eqnarray}
\label{6.3} Y_2 (t) &:=& B(t) - \tfrac{ 1 }{ 2 }  L^{\Gamma} (t)
+ \tfrac{ 1 }{ 2 }  L^{\Delta} (t)  , \qquad  Y_1 (t) :=
Y_2 (t) + \Gamma(t) ,
\nonumber\\[-8pt]\\[-8pt]
   Y_3 (t) &:=& Y_2
(t) - \Delta(t) ,\nonumber
\end{eqnarray}
and note
\begin{equation}
\label{6.5} Y_1 (t)  =  - Y_3 (t)  =
r_1- h t +\tfrac{ 1 }{ 2 }  \bigl( L^{\Gamma} (t) +
L^{\Delta} (t) \bigr)   \ge  \bigl\llvert {Y_2 (t)}\bigr
\rrvert  ,\qquad   0 \le t < \infty .
\end{equation}
We shall show below that, with probability one,
\begin{equation}
\label{6.6} \mbox{\textit{the three-dimensional process} }   \bigl( Y_1
(\cdot), Y_2 (\cdot), Y_3 (\cdot) \bigr)   \mbox{ \textit{exhibits no triple point}}.
\end{equation}
Then the comparisons in (\ref{6.5}) imply
\begin{equation}
\label{6.7} \mathbb{P} \bigl( Y_1 (t) >0 ,  \forall  t \in[0,
\infty) \bigr)   = 1 .
\end{equation}

To prove (\ref{6.6}), it suffices to rule out triple points for the
process $   ( \widehat{Y}_1 (\cdot), \widehat{Y}_2 (\cdot),
\widehat{Y}_3 (\cdot)  ) $ with components $  \widehat{Y}_j
(\cdot) := Y_j (\cdot) + ( L^{\Gamma} (\cdot) - L^{\Delta} (\cdot
) ) /  2 $, $  j=1,2,3, $ namely
\[
\widehat{Y}_1 ( t) =  r_1 - h  t + L^{\Gamma}
(t) \ge  \widehat {Y}_2 ( t) = B(t) \ge  \widehat{Y}_3
( t) =  -r_1 + h  t - L^{\Delta} (t) , \qquad  0 \le t <
\infty .
\]
Consider the set $ E $ of all $  \omega\in\Omega $ for which $
\widehat{Y}_1 ( T) = \widehat{Y}_2 ( T) = \widehat{Y}_3 ( T) =:y $
holds for some $  T = T (\omega) \in(0, \infty) $. Then, for each
$  t \in[0, T) $ we have
\begin{eqnarray*}
y +h  (T- t )  &\ge&  y +h  (T- t ) - L^{\Gamma} (T) + L^{\Gamma}
(t)  =  \widehat{Y}_1 ( t)   \ge  \widehat{Y}_2 ( t)
 =  y+ B(T) - B(t)
\\
&\ge&  \widehat{Y}_3 ( t)  =  y -h  (T- t ) + L^{\Delta}
(T) - L^{\Delta} (t) \ge  y -h  (T- t )
\end{eqnarray*}
%
%
%
thus also
\[
- h   \le  {  B(T) - B(t)   \over T-t}   \le h .
\]

In conjunction with the \textsc{Paley, Wiener and Zygmund} theorem (cf. page
110 in Karatzas and Shreve \cite{MR1121940}), we conclude that $ E $ is included in an
event of $ \mathbb{P}$-measure zero, so (\ref{6.6}) follows. (We are
indebted to \textsc{Dr. Johannes Ruf}, for pointing out the relevance
of the \textsc{Paley--Wiener--Zygmund} theorem here.)

Let us define now
\begin{equation}
\label{6.9} Q_1 (\cdot)  :=  Y_1 (\cdot) ,\qquad
Q_2 (\cdot)  :=  \big| Y_2 (\cdot) \big|
\end{equation}
and apply the companion \textsc{Tanaka} formula (\ref{Tanaka1}) to
get $ Q_2 (\cdot) = \int_0^\cdot\overline{\sgn}  ( Y_2 (t)
) \,\mathrm{d} Y_2 (t) + L^{Q_2} (\cdot) $. In conjunction with
(\ref{6.3}), the fact that $  L^{\Gamma} (\cdot) = L^{Y_1 - Y_2}
(\cdot) $ is flat off the set $  \{ t \ge0 \dvt Y_2 (t) = Y_1 (t) >0\}
 $, and the fact that $  L^{\Delta} (\cdot) = L^{Y_2 - Y_3} (\cdot
) $ is flat off the set $ \{ t \ge0\dvt Y_2 (t) = Y_3 (t) <0\}$, this
leads to %
\begin{equation}
\label{6.10} Q_2 (\cdot)  =  \widehat{V} (\cdot) - \tfrac{ 1 }{ 2 }
  \bigl( L^{Y_1 - Y_2} (\cdot) + L^{Y_1 + Y_2} (\cdot) \bigr)+
L^{Q_2} (\cdot ) ,
\end{equation}
where
\[
\widehat{V} (\cdot)  :=  \int_0^\cdot
\overline{\sgn} \bigl( Y_2 (t) \bigr)\, \mathrm{d} B (t) =  \int
_0^\cdot\sgn \bigl( Y_2 (t) \bigr)
\,\mathrm{d} B (t)
\]
is standard Brownian motion (for this last equality, we have applied
the second property of (\ref{flat}) to the semimartingale $
Y_2(\cdot) $ in (\ref{6.3})). From (\ref{6.9}), (\ref{6.3}) and
(\ref{6.5}), we have then
\begin{eqnarray}
\label{6.100} Q_1 (\cdot) &= & r_1 - h  t + \tfrac{ 1 }{ 2 }   \bigl( L^{Y_1 -
Y_2} (\cdot) + L^{Y_1 + Y_2} (\cdot) \bigr) ,
\\
\label{6.11} 0   &\le&  Q_1 (\cdot)- Q_2 (\cdot) =
 r_1 - h  t - \widehat{V} (\cdot)- L^{Q_2} (\cdot)+ \bigl(
L^{Y_1 - Y_2} (\cdot) + L^{Y_1 +
Y_2} (\cdot) \bigr) .
\end{eqnarray}
But the continuous, increasing process $  L^{Y_1 - Y_2} (\cdot) +
L^{Y_1 + Y_2} (\cdot) $ is flat away from the set
$
\{ Q_1 (\cdot)= Q_2 (\cdot) \} =\{ Y_1 (\cdot)\pm Y_2 (\cdot)=0 \}
 $,
so the theory of the \textsc{Skorokhod} reflection problem gives the
$  \mathbb{P}$-a.e. identities
\[
  L^{Y_1 - Y_2} (\cdot) + L^{Y_1 + Y_2} (\cdot)   \equiv
L^{Q_1 -
Q_2} (\cdot) ,\qquad   \int_0^\infty
\mathbf{ 1}_{\{ Q_1 (t) > Q_2
(t) \} }  \, \mathrm{d} L^{Q_1 - Q_2} (t)  = 0 .
\]
With this identification, the system (\ref{6.11}), (\ref{6.100}) is
written equivalently as
\begin{eqnarray}
Q_1 (t) &=&  r_1 - h  t + \tfrac{ 1 }{ 2 }
L^{Q_1 - Q_2} (t) ,
\\
Q_2 (t) &= & \widehat{V} (t) - \tfrac{ 1 }{ 2 } L^{Q_1 - Q_2}
(t) + L^{Q_2} (t)
\end{eqnarray}
with $  \int_0^\infty\mathbf{ 1}_{\{ Q_2 (t) > 0 \} }  \, \mathrm
{d} L^{ Q_2} (t) =0 $, $  \int_0^\infty\mathbf{ 1}_{\{ Q_1 (t) >
Q_2 (t) \} } \,  \mathrm{d} L^{Q_1 - Q_2} (t) =0  $; that is,
precisely in the form (\ref{R1})--(\ref{R3}) with $  g=0 $, $
r_2=0 $.
Thus, the pairs $  (Q_1 (\cdot), Q_2 (\cdot)) $ and $  (N (\cdot
), M (\cdot)) $ have the same distribution, so the property $
\mathbb{P}
( N (t) >0 ,  \forall  t \in[0, \infty)  )  = 1 $
follows now from (\ref{6.7}) and (\ref{6.9}).


\subsection{\texorpdfstring{Proof of Proposition \protect\ref{Prop6}}{Proof of Proposition 4.4}}
\label{sec62}

In the nondegenerate case $  \sigma^2 <1 $, the properties of
(\ref{ReimWill}) follow from Theorem 1 in
Reiman and Williams \cite{MR921820}; see also Theorem 7.7 in
Bhardwaj and Williams \cite{MR2546421}. In the degenerate case $  \sigma^2 =1 $, they
follow from (\ref{3.4.a}), (\ref{3.18}) -- which give $  A(\cdot) =
\int_0^{  \cdot} \mathbf{ 1}_{ \{ G(t) =0\}}  \, \mathrm{d} A(t) $
and $  \Lambda(\cdot) = \int_0^{  \cdot} \mathbf{ 1}_{ \{ M(t)
=0\}}  \, \mathrm{d} \Lambda(t) $, respectively
-- and the nonattainability of the corner that we just proved.


\section{Questions of uniqueness}
\label{sec7}

Let $  \mathcal{B}(\cdot) = (B_1(\cdot),B_2(\cdot))'  $ be a
planar Brownian motion, and define the matrix- and vector-valued functions
$ \bolds{\Sigma}\dvtx \mathbb{R}^2\to\mathbb{R}^{2 \times2}  $ and
$  \bolds{\mu}
\dvtx \mathbb{R}^2\to\mathbb{R}^2  $ by
\begin{eqnarray*}
\bolds{\Sigma}(x) &=&  \mathbf{ 1}_{\{x_1> x_2\}}  \operatorname{diag}(\rho ,
\sigma)+\mathbf{ 1}_{ \{x_1\leq x_2\}}   \operatorname{diag}(\sigma ,\rho),
\\
\bolds{\mu}(x) &=&   \mathbf{ 1}_{\{x_1>x_2\}}   (-h,g)'+
\mathbf{ 1}_{ \{
x_1\leq x_2\}}  (g,-h)'
\end{eqnarray*}
for $  x = (x_1, x_2) \in\mathbb{R}^2 $. We are interested in
questions of uniqueness for the system of stochastic differential
equations (\ref{1.3.a}) and (\ref{1.4.a}), written now a bit more
conveniently in the vector form
\begin{equation}
\label{10.1} \mathrm{d}\mathcal{X}(t) = \bolds{\mu} \bigl(\mathcal{X}(t)
\bigr)\,  \mathrm{d}t+\bolds{\Sigma} \bigl(\mathcal{X}(t) \bigr) \, \mathrm {d}
\mathcal {B}(t)+\mathrm{d}\mathcal{L}^{\mathcal{X}}(t) .
\end{equation}
Here we have denoted the vector of semimartingale local time processes
accumulated at the origin by the components of the planar process $
\mathcal{X}(\cdot) =  ( X_1 (\cdot), X_2 (\cdot)  )' $ as
\begin{equation}
\label{10.5} \mathcal{L}^{\mathcal{X}}(\cdot) :=  \bigl(L^{X_1} (
\cdot),L^{X_2} (\cdot) \bigr)' ;
\end{equation}
these local times are responsible for keeping the planar process $
\mathcal{X}(\cdot)  $ in the nonnegative quadrant.

\subsection{Pathwise uniqueness}

For any given initial condition $ \mathcal{X}(0) $ in the punctured
nonnegative quadrant of (\ref{1.2}), we have shown that the stochastic
differential equation \eqref{10.1} has a solution. We want to show
that this solution is pathwise unique, up to the first hitting time of
the corner of the quadrant. Under the condition (\ref{SigmaGeRho}),
the origin is never hit by the process $  \mathcal{X}(\cdot) $, so
pathwise uniqueness will hold then for all times. The key step is to
define the new planar process
\begin{equation}
\label{10.2} \mathcal{Z}(\cdot) = \bigl( Z_1 (\cdot),
Z_2 (\cdot) \bigr)'   \qquad     \mbox{with }
     Z_i (\cdot) := X_i (\cdot) - L^{ X_i } (
\cdot)  \qquad    (i=1, 2) ,
\end{equation}
and note, from the \textsc{Skorokhod} reflection problem once again, that
\begin{equation}
\label{10.3} L^{ X_i }( t)  =  Z^*_i (t) :=
\max_{0 \le s \le t} \bigl( - Z_i (s) \bigr) \vee0 ,\qquad   0 \le t <
\infty .
\end{equation}
(In particular, each $  X_i (\cdot) $, $  i=1, 2 $ is the \textsc
{Skorokhod} reflection of the semimartingale $  Z_i (\cdot) $ in
(\ref{10.2}).) We observe that $  \mathcal{X} (\cdot) $ \textit{solves the equation \textup{(\ref{10.1})}, if and only if $  \mathcal{Z}
(\cdot) $ solves the stochastic differential equation with
path-dependent coefficients}
\begin{equation}
\label{10.4} \mathrm{d}\mathcal{Z}(t) = \bolds{\mu} \bigl(\mathcal{Z}(t) +
\mathcal {Z}^*(t) \bigr) \,\mathrm{d}t+\bolds{\Sigma} \bigl(\mathcal{Z}(t) +
\mathcal {Z}^*(t) \bigr) \, \mathrm{d}\mathcal{B}(t) ,\qquad   \mathcal{Z}(0) =
\mathcal{X}(0) \in\mathfrak{ S},
\end{equation}
where, with the notation of (\ref{10.3}), we have set
\begin{equation}
\label{10.6} \mathcal{Z}^*(t)  :=  \bigl(  Z_1^* (t),
Z_2^* (t) \bigr)' =  \Bigl( \max_{0 \le s \le t}
\bigl( - Z_1 (t) \bigr) \vee0 ,   \max_{0 \le s \le t} \bigl( -
Z_2 (t) \bigr) \vee0 \Bigr)' .
\end{equation}

Indeed, if $  (\mathcal{X}(\cdot), \mathcal{B}(\cdot)) $ is a
solution of (\ref{10.1}) with the vector process $  L^{\mathcal{X}}
(\cdot) $ given as in (\ref{10.5}), and if we define $  \mathcal
{Z}(\cdot)  $ in accordance with (\ref{10.2}), then $  (\mathcal
{Z}(\cdot), \mathcal{B}(\cdot)) $ is a solution of (\ref{10.4}).
And conversely, if $  (\mathcal{Z}(\cdot), \mathcal{B}(\cdot)) $
is a solution of (\ref{10.4}) and we define
\begin{equation}
\label{10.7} \mathcal{X} (\cdot)  :=  \mathcal{Z} (\cdot) + \mathcal{Z}^* (
\cdot)
\end{equation}
in the notation of (\ref{10.6}), we identify $  \mathcal{Z}^* (\cdot
) $ as the coordinate-wise local time of $  \mathcal{X} (\cdot)
$, namely $  \mathcal{Z}^* (\cdot) = \mathcal{L}^{\mathcal{X}}
(\cdot) $ as in (\ref{10.5}). In particular, pathwise uniqueness for
the equation
\eqref{10.4} with path-dependent coefficients, implies pathwise uniqueness
for the equation \eqref{10.1} with local times.

For the equation \eqref{10.4}, this pathwise uniqueness result can be
seen by exploiting the fact that the only critical time-points occur
when the process $  \mathcal{X} (\cdot) $ of (\ref{10.7}) hits
either the boundary of the quadrant, or its diagonal. Thus, we set $
\tau_0 :=0 $, and introduce inductively the stopping times
\begin{eqnarray*}
\tau_{2n+1} &:=&  \inf \bigl\{  t > \tau_{2n} \dvt  \bigl(
Z_1 (t)+ Z^*_1(t) \bigr) \bigl( Z_2 (t)+
Z^*_2(t) \bigr) = 0   \bigr\},
\\
\tau_{2n+2} &:=&  \inf \bigl\{  t > \tau_{2n+1} \dvt
Z_1 (t)+ Z^*_1(t) = Z_2 (t)+
Z^*_2(t)   \bigr\}
\end{eqnarray*}
for $  n \in\mathbb{N}_0 $. It might happen that $\tau_1=\tau_0=0$ holds with positive probability, if $ \mathcal{X}(0)=\mathcal
{Z}(0) $ is on one of the faces of the quadrant. Apart from that, we
have $ \tau_{k+1}>\tau_k  $ for all $ k\geq1 $; and with $
\tau:=\sup_{k \in\mathbb{N}}  \tau_k $, we observe that $
\mathcal{X}(\tau)=\mathcal{Z}(\tau)+\mathcal{Z}^*(\tau)$ lies on
one of the faces of the quadrant and also on its diagonal, hence $
\mathcal{X}(\tau)=(0,0) $. Whereas, under the condition (\ref
{SigmaGeRho}), we know that $ \mathcal{X} (\cdot) $ never reaches
the corner of the quadrant, so $ \tau  $ is almost surely infinite.
Then pathwise uniqueness
up to $ \tau $ can be established by induction; namely, by showing
that on each time-interval $ [\tau_{k-1},\tau_k] $ with $  k \ge
1 $ the process $  \mathcal{Z} (\cdot) $ solves an equation with
pathwise unique solution.

Assume first that $ k-1 $ is odd; then on the time-interval $ [\tau_{k-1},\tau_k] $ the drift and
diffusion coefficients are constant, so pathwise uniqueness follows
immediately. If $  k-1  $ is even, then $  \mathcal{Z}^* (\cdot)
$ is
constant on the interval $ [\tau_{k-1},\tau_k] $, and the process
$
\mathcal{X} ( \cdot) =\mathcal{Z} ( \cdot) + \mathcal{Z}^* ( \cdot
) $
solves on $  [\tau_{k-1},\tau_k] $ the equation studied in
Fernholz \textit{et al.} \cite{2011arXiv1108.3992F}, Theorem 5.1, where it was shown that the
solution of this equation is pathwise unique. This completes the
induction argument.

Invoking the \textsc{Yamada--Watanabe} theory (e.g.,
Karatzas and Shreve \cite{MR1121940}, pages 308--311), we obtain from all this the following result.

\begin{proposition}
\label{Prop4}
For the system of equations (\ref{1.3.a}) and (\ref{1.4.a}), pathwise
uniqueness holds up until the first hitting time $  \tau $ of (\ref{1.5}).

Under the condition $  1/2 \le\sigma^2 <1 $ of (\ref{SigmaGeRho})
pathwise uniqueness, thus also uniqueness in distribution, hold for all
times, so the solution of (\ref{1.3.a}) and (\ref{1.4.a}) is then
strong; that is, for all $ 0 \le t < \infty $ we have
\begin{equation}
\label{Incl1} \mathfrak{F}^{ (X_1, X_2)} (t)   \subseteq
\mathfrak{F}^{ (B_1,
B_2)} (t) .
\end{equation}
\end{proposition}

This completes the proof of Theorems \ref{Thm1} and \ref{Thm2}. It is
an open question to settle, whether such pathwise uniqueness and
strength hold also past the first hitting time of the corner, or not.


\subsection{Uniqueness in distribution}
\label{subsec8}

This section will be devoted to the proof of Theorem \ref{Thm3} on
uniqueness in distribution. In the light of Proposition \ref{Prop4}
and Theorem \ref{Thm2}, this needs elaboration only when $  0<\sigma^2 < 1 / 2 $ as in
(\ref{RhoGreSigma}). In this case, the state
process $  \mathcal{X}(\cdot) =  ( X_1 (\cdot), X_2 (\cdot)
)'  $ of (\ref{10.1}) can reach the corner of the quadrant $
[0, \infty)^2 $ in finite time, and the question is whether it can be
continued beyond that time in a well-defined and unique-in-distribution manner.

\begin{pf*}{Proof of Theorem \ref{Thm3}}
 It is quite straightforward to see
that the weak solution construction of Section \ref{sec4}, culminating
with the continuous, nonnegative processes $  X_1 (\cdot),  X_2
(\cdot)  $ defined in Section \ref{sec4.4}, makes perfectly good
sense also for an initial condition $  (X_1 (0), X_2 (0)) = (x_1, x_2)
=(0,0 ) $ at the corner of the quadrant. Together with the independent
Brownian motions $  B_1 (\cdot),  B_2 (\cdot)  $ of (\ref{3.7}),
(\ref{3.7.a}), these processes $  X_1 (\cdot),  X_2 (\cdot)  $
are constituents of a weak solution for the system (\ref{10.1}), and
from (\ref{10.5})--(\ref{10.3}) we have the representations
\[
X_i (\cdot)  =  Z_i (\cdot)+ \max_{0 \le s \le  \cdot}
\bigl( - Z_i (s) \bigr)  , \qquad   i=1,2  .
\]
Now, on an extension of the filtered probability space on which $  (
B_1 (\cdot),  B_2 (\cdot) )  $ is still planar Brownian motion, we
unfold these continuous, nonnegative semimartingales as
\[
X_i(\cdot) = \big\vert Z_i^\flat(\cdot) \big\vert
,\qquad   \mbox{where }   Z_i^\flat(\cdot) = \int
^{ \cdot}_0 \overline{\operatorname{sgn}}
\bigl(Z_i^\flat(s) \bigr)  \, \mathrm d Z_i(s)
.
\]
This is done using the %
Prokaj \cite{MR2494646} construction of Section \ref{sec4.1.a} once
again. It follows that the process $  \mathcal Z^\flat(\cdot) := (
Z_1^\flat(\cdot), Z_2^\flat(\cdot))'  $ with component-wise
absolute values $  \vert\mathcal Z^\flat(\cdot) \vert: =
(\vert Z_1^\flat(\cdot)\vert  , \vert Z_2^\flat(\cdot) \vert
)'$ $ = \mathcal X(\cdot)   $ satisfies the vector stochastic
differential equation
\begin{eqnarray}
\label{eq: unfolded} %
  \mathrm d \mathcal
Z^\flat(\cdot) &   = & \mathbf{I}\bigl(\mathcal Z^\flat (t)\bigr)
  \bigl( \bolds{\mu}\bigl(\mathcal X(t) \bigr)\,  \mathrm {d} t + \bolds{\Sigma }\bigl(
\mathcal X(t)\bigr)  \,\mathrm {d} \mathcal{B}(t) \bigr)
\nonumber\\[-8pt]\\[-8pt]
&   = & \mathbf{I}\bigl(\mathcal Z^\flat(t)\bigr)   \bolds{\mu}\bigl(\big\vert
\mathcal {Z}^\flat(t) \big\vert\bigr)  \,\mathrm {d} t + \mathbf{I}\bigl(\mathcal
Z^\flat(t)\bigr)   \bolds{\Sigma}\bigl(\big\vert\mathcal{Z}^\flat(t)
\big\vert\bigr) \,  \mathrm {d} \mathcal{B}(t)   ,\qquad   0 \le t < \infty ,\nonumber
\end{eqnarray}
with the initial condition $  \mathcal Z^\flat(0) = 0  $, the
indicator matrix function
\[
  \mathbf{I}(z) := \operatorname{diag} \bigl( \overline{\operatorname{sgn}} (z_1)
, \overline{\operatorname{sgn}} (z_2) \bigr) \qquad   \mbox{and the notation }
     z =(z_1, z_2)' \in
\mathbb{R}^2  ,    |z| := \bigl( |z_1|, |z_2|\bigr)'
 .
\]
The functions $  \mathbf{I}(\cdot)  \bolds{\mu}(\vert  \cdot
\vert)  $ and $  \mathbf{I}(\cdot)  \bolds{\Sigma}(\vert  \cdot
\vert)  $ are piecewise constant in the interior of each one of the
eight wedges
\begin{eqnarray*}
&&\bigl\{(z_1, z_2)  \dvt   z_1 \ge0   ,
z_2 \ge0  , z_1 > z_2\bigr\}   ,\qquad
\bigl\{(z_1, z_2)  \dvt   z_1 \ge0   ,
z_2 \ge0  , z_1 \le z_2\bigr\}   ,
\\
&&\bigl\{(z_1, z_2)  \dvt   z_1 < 0   ,
z_2 \ge0  , z_1 \ge z_2\bigr\}   ,\qquad
\bigl\{(z_1, z_2)  \dvt   z_1 < 0   ,
z_2 \ge0  , z_1 \le z_2\bigr\}   ,
\\
&&\bigl\{(z_1, z_2)  \dvt   z_1 \ge0   ,
z_2 < 0  , z_1 > z_2\bigr\}   ,\qquad
\bigl\{(z_1, z_2)  \dvt   z_1 \ge0   ,
z_2 < 0  , z_1 \le z_2\bigr\}   ,
\\
&&\bigl\{(z_1, z_2)  \dvt   z_1 < 0   ,
z_2 < 0  , z_1 > z_2\bigr\}   ,\qquad
\bigl\{(z_1, z_2)  \dvt   z_1 < 0   ,
z_2 < 0  , z_1 \le z_2\bigr\}   .
\end{eqnarray*}
Theorem 2.1 of Bass and Pardoux \cite{MR917679} (see also Theorem 5.5 in
Krylov \cite{MR2078536}, as well as Exercise~7.3.4, pages 193--194 in
Stroock and Varadhan \cite{MR532498}) guarantees that uniqueness in distribution holds for the
system of equations (\ref{eq: unfolded}) with piecewise constant
coefficients. (For the applicability of this result, the nondegeneracy
condition $\rho  \sigma>0 $ is crucial.)\vadjust{\goodbreak} Whereas, because the
distribution of $  \mathcal Z^\flat(\cdot)   $ is uniquely determined
from (\ref{eq: unfolded}), it is checked fairly easily that the distribution
of $  \mathcal X(\cdot) = \vert\mathcal Z^\flat(\cdot) \vert
$ is
uniquely determined from~(\ref{10.1}).
\end{pf*}

The proofs of all three Theorems \ref{Thm1}--\ref{Thm2} are now
complete. We have
shown in particular that, under the condition $  0<\sigma^2 < 1 / 2
$ as in
(\ref{RhoGreSigma}), the planar process $  \mathcal X(\cdot) $ can hit
the corner of the quadrant but then ``finds a way to extricate itself''
in such a manner that uniqueness in distribution holds. This aspect of the
diffusion is reminiscent of Section 3 of %
Bass and Pardoux \cite{MR917679}; we will see in the \hyperref[app]{Appendix} that these features hold
also in
the other degenerate case $  \sigma=0 $, $  \rho=1 $.

%
\section{Alternative systems and filtration identities}
\label{sec9}

Throughout this section, we shall place ourselves under the condition
\[
1/2   \le  \sigma^2   <  1
\]
for simplicity. We disentangle the pair $  (B_1 (\cdot), B_2(\cdot
)) $ from $  (W_1 (\cdot), W_2(\cdot)) $ in (\ref{2.3}) and (\ref
{2.4}), and rewrite (\ref{3.11}), (\ref{3.12}) in the form of a
system of equations driven by the planar Brownian motion $  \mathcal
{W} (\cdot) =(W_1(\cdot),   W_2 (\cdot))' $, namely
\begin{eqnarray}\label{3.11.a}
X_1 (t) &=& x_1 + \int_0^t
( g   \mathbf{ 1}_{ \{ X_1 (s) \le
X_2 (s) \} } - h   \mathbf{ 1}_{ \{ X_1 (s) > X_2 (s) \} } )
\,\mathrm{d} s +   \rho\int_0^t \mathbf{
1}_{ \{ X_1 (s) > X_2 (s)
\} }  \,\mathrm{d} W_1 (s)\qquad
\nonumber\\[-8pt]\\[-8pt]
    &&{}             +  \int_0^t
\mathbf{ 1}_{ \{ X_1 (s) \le X_2 (s) \}
}  \bigl( \sigma \, \mathrm{d} W_2 (s) +
\mathrm{d} L^{X_1 \wedge
X_2} (s) \bigr) ,\qquad  0 \le t < \infty,\nonumber
\\\label{3.12.a}
X_2 (t) &=&  x_2 + \int_0^t
( g   \mathbf{ 1}_{ \{ X_1 (s) > X_2 (s) \} } - h   \mathbf{ 1}_{ \{ X_1 (s) \le X_2 (s) \} } )
\,\mathrm{d} s  -   \rho \int_0^t \mathbf{
1}_{ \{ X_1 (s) \le X_2 (s) \} } \, \mathrm{d} W_1 (s)
\nonumber\\[-8pt]\\[-8pt]
   &&{}             +  \int_0^t
\mathbf{ 1}_{ \{ X_1 (s) > X_2 (s) \}
}  \bigl( - \sigma \, \mathrm{d} W_2 (s) +
\mathrm{d} L^{X_1 \wedge
X_2} (s) \bigr) ,\qquad  0 \le t < \infty .\nonumber
\end{eqnarray}

Repeating the argument of Proposition \ref{Prop4} and using once again
pathwise uniqueness results from Theorem 4.2 in
Fernholz \textit{et al.} \cite{2011arXiv1108.3992F}, one can show the \textit{pathwise uniqueness
and strong solvability of the system \textup{(\ref{3.11.a})}  and  \textup{(\ref{3.12.a})},
under the condition \textup{(\ref{SigmaGeRho})}.} In particular, we have
\begin{equation}
\label{Incl2} \mathfrak{F}^{ (X_1, X_2)} (t)   \subseteq
\mathfrak{F}^{ (W_1,
W_2)} (t)  \qquad   \forall  0 \le t < \infty  .
\end{equation}

\noindent$\bullet $ In a similar manner, recalling the expressions for $  (B_1
(\cdot), B_2(\cdot)) $ in terms of $  (V_1 (\cdot), V_2(\cdot))
$ in (\ref{3.7}), (\ref{3.7.a}) and disentangling the former from the
latter, %
we can rewrite the system of equations (\ref{3.11}), (\ref{3.12}) in
the form
\begin{eqnarray}\label{3.11.b}
X_1 (t) &= &x_1 + \int_0^t
( g   \mathbf{ 1}_{ \{ X_1 (s) \le
X_2 (s) \} } - h   \mathbf{ 1}_{ \{ X_1 (s) > X_2 (s) \} } )
\,\mathrm{d} s +   \rho\int_0^t \mathbf{
1}_{ \{ X_1 (s) > X_2 (s)
\} } \, \mathrm{d} V_1 (s)\qquad
\nonumber\\[-8pt]\\[-8pt]
  &&{}              +  \int_0^t
\mathbf{ 1}_{ \{ X_1 (s) \le X_2 (s) \}
}  \bigl( \sigma\,  \mathrm{d} V_2 (s) +
\mathrm{d} L^{X_1 \wedge
X_2} (s) \bigr) ,\qquad  0 \le t < \infty,\nonumber
\\\label{3.12.b}
X_2 (t) &=&  x_2 + \int_0^t
( g   \mathbf{ 1}_{ \{ X_1 (s) > X_2 (s) \} } - h   \mathbf{ 1}_{ \{ X_1 (s) \le X_2 (s) \} } )
\,\mathrm{d} s  +   \rho \int_0^t \mathbf{
1}_{ \{ X_1 (s) \le X_2 (s) \} } \, \mathrm{d} V_1 (s)
\nonumber\\[-8pt]\\[-8pt]
  &&{}              +  \int_0^t
\mathbf{ 1}_{ \{ X_1 (s) > X_2 (s) \} }   \bigl( \sigma \, \mathrm{d} V_2 (s) +
\mathrm{d} L^{X_1 \wedge
X_2} (s) \bigr) ,\qquad  0 \le t < \infty .\nonumber
\end{eqnarray}
\textit{This system admits a unique-in-distribution weak solution} (recall
Remark \ref{31}), \textit{but not a strong one.} Indeed, pathwise uniqueness
cannot hold for (\ref{3.11.b}) and (\ref{3.12.b}) if the process $
\mathcal{X} (\cdot) $ hits the diagonal of the quadrant; but the
diagonal is hit with positive probability during any time-interval $
(0,t) $ with $ 0 < t < \infty $, so in conjunction with Remark \ref
{51} we have (with strict inclusion)
\begin{equation}
\label{Incl3} \mathfrak{F}^{ (V_1, V_2)} (t)  =  \mathfrak{F}^{ (R_1, R_2)}
(t)   \subsetneqq  \mathfrak{F}^{
(X_1, X_2)} (t)  .
\end{equation}

\subsection{Filtration identities}
\label{sec9a}

Let us recall the equations of (\ref{3.7}) and (\ref{3.7.a}), written
now a
bit more conspicuously in the form
\begin{eqnarray*}
B_1(t)  &=&  \int_0^t \mathbf{
1}_{ \{ X_1 (s) > X_2 (s)\} }  \, \mathrm{d} V_1 (s) + \int_0^t
\mathbf{ 1}_{ \{ X_1 (s) \le X_2 (s)
\} }  \, \mathrm{d} V_2 (s)  ,
\\
B_2(t)  &=&   \int_0^t \mathbf{
1}_{ \{ X_1 (s) \le X_2 (s) \} }  \, \mathrm{d} V_1 (s) + \int_0^t
\mathbf{ 1}_{ \{ X_1 (s) > X_2 (s) \}
}  \, \mathrm{d} V_2 (s)  .
\end{eqnarray*}
From these equations and the filtration inclusion $  \mathfrak{F}^{
(V_1, V_2)} (t) \subseteq\mathfrak{F}^{ (X_1, X_2)} (t)  $ of
(\ref{Incl3}), we conclude that the reverse inclusion $  \mathfrak
{F}^{ (B_1, B_2)} (t) \subseteq\mathfrak{F}^{ (X_1, X_2)} (t)  $
of (\ref{Incl1}) also holds. We have thus argued, for all $  0 \le t
< \infty $, the filtration identity
\begin{equation}
\label{Incl4} \mathfrak{F}^{ (X_1, X_2)} (t)   =   \mathfrak{F}^{ (B_1, B_2)}
(t)  .
\end{equation}

Similar reasoning, applied to the equations
\[
  W_1 (t) = \int_0^t
\operatorname{sgn} \bigl( X_1(s) - X_2 (s) \bigr)\,   \mathrm{d}
V_1 (s) \quad  \mathrm{and} \quad   W_2 (t) = - \int
_0^t \operatorname{sgn}   \bigl( X_1(s)
- X_2 (s) \bigr)  \, \mathrm{d} V_2 (s)
\]
from (\ref{3.8}), (\ref{3.8.a}), leads in conjunction with (\ref
{Incl2}) to the filtration identity
\begin{equation}
\label{Incl5} \mathfrak{F}^{ (X_1, X_2)} (t)   =   \mathfrak{F}^{ (W_1, W_2)}
(t)  , \qquad  0 \le t < \infty .
\end{equation}

\noindent$\bullet $
We conclude from these two identities and (\ref{Incl3}) that, \textit{under the condition \textup{(\ref{SigmaGeRho})}, both of the planar Brownian
motions $  (B_1(\cdot),B_2(\cdot)) $ and $  (W_1(\cdot),
W_2(\cdot)) $ generate the same filtration as the state-process} $
(X_1(\cdot), X_2(\cdot)) $; whereas the planar Brownian motion $
(V_1(\cdot),V_2(\cdot)) $ generates a \textit{strictly smaller}
filtration, that of the ranked processes $  (R_1(\cdot),R_2(\cdot))
$ in (\ref{ranks}).

\section{Questions of recurrence and transience}
\label{sec10}

Hobson and Rogers \cite{MR1198420} (see also Dupuis and Williams \cite{MR1288127} and Chen \cite{MR1410114})
study a
reflecting Brownian motion $  Z(t) := (X (t), Y(t)) $, $  0 \le t <
\infty $ where the co\"ordinate processes satisfy the equations
\begin{eqnarray}
\label{eq: hob 1} X(t)  &=&  x + B(t) + \mu t + \alpha L^Y (t)+
L^X (t)  ,
\\
\label{eq: hob 2} Y(t)  &=&  y + W(t) + \nu t + \beta L^X (t)+
L^Y (t)   .
\end{eqnarray}
Here $  \alpha  ,   \beta  ,   \mu ,   \nu  $ are fixed
real numbers with $  (\mu, \nu) \neq(0,0) $; the process $  (B
(\cdot) , W (\cdot)) $ is a planar Brownian motion with nonsingular
covariance; and $  L^X (\cdot) $ (resp., $  L^Y (\cdot)  $)
is the local time process at the origin of $  X (\cdot) $ (resp., $  Y (\cdot) $).

Consider a bounded neighborhood $  \mathcal N $ of the origin in $  [0,
\infty)^2 $, and let $  \mathbf{T}:=   \inf\{ s \ge0\dvt (X(s), Y(s)) \in
{\mathcal N} \} $ be the first entry time in $  \mathcal N $.
Theorem 1.1
in Hobson and Rogers \cite{MR1198420} provides a classification of recurrence and
transience for the reflecting Brownian motion:
\begin{enumerate}[2.]
\item[1.] For every initial point $  Z(0) = z \in[0, \infty)^2 $, we have
\[
\mu+ \alpha  \nu^{-} \le0  ,\qquad   \nu+ \beta  \mu^- \le0
\quad \Longleftrightarrow \quad   \mathbb{P}^z (\mathbf{T} < \infty) = 1   .
\]
\item[2.]\label{HB:it:2} For every initial point $  Z(0) = z \in[0,
\infty)^2  $ and for some constant $  C > 0 $, we have
\[
\mathbb{E}^z (\mathbf{T}) < \infty \quad  \Longleftrightarrow \quad  \mu+
\alpha  \nu^{-} <   0, \qquad  \nu+ \beta  \mu^- < 0
\quad \Longrightarrow \quad   \mathbb{E}^z(\mathbf{T}) \le C\bigl(1 + \Vert z \Vert\bigr)
  .
\]
\end{enumerate}
Here the dichotomies are determined by the \textit{effective drift rates}
$  \mu
+ \alpha  \nu^- $ and $  \nu+ \beta  \mu^- $, rather than the pure
drift rates $  \mu $ and $  \nu $. This is because the local time
$  L^X
(t) $ grows like $  \mu^-   t  $, so the effective drift rate of the
process $  Y(\cdot) $ is $  \nu+ \beta\mu^- $; similarly, the local
time $  L^Y (t) $ grows like $  \nu^-  t   $, so the effective
drift rate
of the process $  X (\cdot) $ is $  \mu+ \alpha\nu^- $.

\noindent$\bullet $
Let us apply this result to the system of (\ref{3.16}) and (\ref{3.17})
for the ranks $  R_1 (\cdot) = X_1 (\cdot) \vee X_2 (\cdot)  $,
$  R_2 (\cdot)= X_1 (\cdot) \wedge X_2 (\cdot) $ of the processes
$  X_1(\cdot),   X_2 (\cdot) $ constructed in (\ref{3.11}) and (\ref
{3.12}), assuming
\[
  \lambda >  0 \quad  \mbox{and} \quad  0   < \sigma < 1
\]
in (\ref{1.1}). Comparing (\ref{eq: hob 1})--(\ref{eq: hob 2}) with
\begin{eqnarray}
\label{HR} R_1 (t)- R_2 (t)  &=&  ( r_1 -
r_2) - \lambda  t + \rho  V_1(t) - \sigma
V_2(t) + L^{R_1 - R_2} (t) - L^{ R_2} (t) ,
\\
\label{2.23.a} R_2 (t) &=&  r_2 + g  t +
\sigma  V_2(t)- \bigl( L^{R_1 - R_2} (t) / 2 \bigr) +
L^{R_2} (t) ,
\end{eqnarray}
an equivalent form in which the system of (\ref{2.20})--(\ref{2.21})
can be cast, we make the identifications
\[
  \mu= - \lambda= -(g+h) < 0  , \qquad  \nu= g > 0 , \qquad   \alpha= -1 ,\qquad
 \beta= ( - 1/ 2 ) .
\]
The first passage time takes the form $  \mathbf{T}:= \inf\{s \ge0\dvt
(R_1(s) - R_2(s), R_2(s) ) \in{\mathcal N}  \} $, and the effective
drift rates become
\[
- (g + h) + (-1) \cdot0 = - (g+h) < 0  , \qquad   g + ( - 1/ 2 )   (g + h) = (g -
h) /2   .
\]
Thus, from %
Hobson and Rogers \cite{MR1198420}, for every initial point $  (r_1 - r_2, r_2) \in(0,
\infty)^2 $ we have:
\begin{enumerate}[]
\item[]$  \mathbb{P}( \mathbf{T} < \infty) < 1 $, if $  g > h $;

\item[]$  \mathbb{P}( \mathbf{T} < \infty) = 1 $, if $  g \le h $; and

\item[]$  \mathbb{E}(\mathbf{T}) < \infty $, if $  g < h   $.
\end{enumerate}
We translate this observation to the following claims for $
(X_1(\cdot), X_2(\cdot)) $:
\begin{itemize}[$\bullet$]
\item[$\bullet$]\textit{If $  g > h $, we have $  \mathbb{P}   ( \inf_{0
\le t <
\infty} R_1(t) = 0 ) < 1 $};
\item[$\bullet$]\textit{If $  g \le h $, we have $  \mathbb{P}   ( \inf_{0
\le t
< \infty} R_1(t) = 0 ) = 1 $}.
\textit{Furthermore, the condition $ g<h $ is necessary and sufficient
for positive recurrence, that is, for the hitting time of any given
Borel set with positive Lebesgue measure by the vector process $
(X_1(\cdot), X_2(\cdot)) $ to have finite expectation for all
starting points $ z\in[0,\infty)^2 $.}
\end{itemize}
The reflection and covariance matrices in (\ref{HR})--(\ref{2.23.a})
do not
satisfy in general the so-called ``skew-symmetry condition'' of
Harrison and Williams \cite{MR877593}, Williams \cite{MR894900}. Hence, the unique invariant distribution
is not of ``exponential form'' in general. It is an open problem to identify
the general form of the invariant distribution for the process $
(R_1(\cdot) -
R_2(\cdot), R_2(\cdot)) $; but we describe this invariant
distribution in a
special case in the subsection that follows, based on results of
Dieker and Moriarty \cite{MR2472171}.

%
\subsection{Densities of sum-of-exponentials type}
\label{sec10.1}

\begin{figure}

\includegraphics{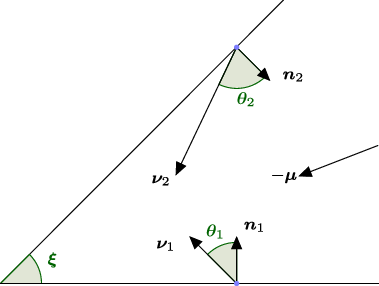}

\caption{The wedge with the reflection parameters and the drift.}\vspace*{-2pt}
\label{fig:2}
\end{figure}

In this subsection, we shall study a special type of invariant
densities for
the ranks $  (R_1(\cdot), R_2(\cdot))  $, applying Theorem 1 of
Dieker and Moriarty \cite{MR2472171}.
This result provides, in certain cases, a formula for the invariant
probability density of a reflected Brownian motion with drift, in a
wedge with
oblique constant reflection on its faces. If
$  \bolds{\alpha}=( \bolds{\theta}_1+\bolds{\theta}_2) / \bolds{\xi}=-\ell $ for
some integer
$ \ell\geq0 $,
then the invariant
probability density $ \mathfrak{ p} (\cdot , \cdot) $ is of the
sum-of-exponentials type, that is, proportional to
\begin{equation}
\label{eq: sum-of-exponential 1} \pi(x) = \sum_{k=0}^\ell
c_k \bigl[{ \bigl\langle\bolds{\mu},{(\mathbf{I}_2 - \mathrm{Rot}_k)\bolds{\nu}_1 } \bigr\rangle  \mathrm{e}^{ -
\langle\bolds{\mu},{(\mathbf{I}_2 - \mathrm{Rot}_k ) x} \rangle} -
\bigl\langle\bolds{\mu},{(\mathbf{I}_2 - \mathrm{Ref}_k)  \bolds{\nu
}_1} \bigr\rangle   \mathrm{e}^{ -  \langle{\bolds{\mu}},{(\mathbf{I}_2 -
\mathrm{Ref}_k ) x} \rangle}} \bigr]
\end{equation}
for $  x \in\{ (x_1, x_2) \in(0, \infty)^2  :   0 < x_2 < x_1
\tan(\bolds{\xi})\}  $. Here $  \mathrm{Rot}_k  $ and $  \mathrm{Ref}_k  $ are rotation and reflection matrices, respectively,
\begin{eqnarray*}
\mathrm{Rot}_k &: = &- %
\pmatrix{ \cos(2\bolds{\theta}_1 + 2 k \bolds{\xi}) & - \sin(2\bolds{\theta}_1 + 2 k {
\bolds{\xi}})
\cr
\sin(2\bolds{\theta}_1 + 2 k \bolds{\xi}) & \cos(2\bolds{\theta}_1
+ 2 k \bolds{\xi})} %
  , \qquad   k = 0, 1, \ldots, \ell  ,
\\
\mathrm{Ref}_k&:=& \mathrm{Rot}_k   \mbox{J} ,
\\
\mathrm{J}&:=&\operatorname{diag}(1,-1)
\end{eqnarray*}
and
\[
 c_k   :=   (-1)^k \cdot
\frac{\prod_{0 \le i < j \le
\ell,
i,j \neq k}  \langle\bolds{\mu},{(\mathrm{Rot}_i - \mathrm{Rot}_j)
\mathbf{e}_1} \rangle} {
\langle\bolds{\mu},{(\mathrm{Ref}_k - \mathrm{Rot}_k)\bolds{\nu}_1
} \rangle}  , \qquad   k = 0, 1, \ldots, \ell  .
\]
In formula \eqref{eq: sum-of-exponential 1} and in the definition of
$ c_k $, the
vector $ \bolds{\nu}_i $ denotes the reflection vector on the $
i $th
face; see
Figure \ref{fig:2}. This vector is
usually normalized so that $ \langle{\bolds{\nu}_i},{\mathbf{n}_i} \rangle=1$. Note, however,
that this
normalization is made just for convenience; it does not affect either
the reflected
process itself, or even the formula \eqref{eq: sum-of-exponential 1}
(as its effect cancels out). Thus, we can safely apply the result with
unnormalized reflection vectors that we have already computed in
Section \ref{sec6}.

To apply Theorem 1 of
Dieker and Moriarty \cite{MR2472171}, first we transform
\eqref{2.16} and \eqref{2.17} into $  \overline{R}_1(\cdot)=R_1(\cdot
)/\rho $,
$  \overline{R}_2(\cdot)=R_2(\cdot)/\sigma $. The resulting process
$ ( \overline{R}_1(\cdot), \overline{R}_2(\cdot)) $ takes values
in 
\[
\bigl\{{(x,y)\in\mathbb{R}^2 \dvt 0\leq y\leq\tan(\bolds{\xi})   x} \bigr
\},
\]
where we use the notation $ \cos(\bolds{\xi})=\sigma $ from \eqref
{eq: XI}. The data of the reflection problem become\vspace*{-1pt} %
\begin{eqnarray*}
\bolds{\mu}&=&(h/\rho,-g/\sigma)',\qquad  \mathbf{n}_1=\bolds{\nu}_1 =(0,1)',\qquad  \mathbf{n}_2 =( \rho, -
\sigma)',\qquad  \bolds{\nu}_2=(1/\rho,-1/\sigma)',
\\
{\bolds{\theta}}_1&=& 0, \qquad {\bolds{\theta}}_2=2  \bolds{\xi}-
\frac{ \uppi }{2}, \qquad  \bolds{\alpha}=\frac{ {\bolds{\theta}}_1+{\bolds{\theta}}_2 }{\bolds{\xi}}=2-\frac
{\uppi}{ 2 \bolds{\xi} }.
\end{eqnarray*}
The result of \textsc{Dieker and Moriarty} can be applied if $ \bm
\alpha $ is a
nonpositive integer $ -\ell $, and this amounts to %
$  \bolds{\xi}=\frac{\uppi}{2(\ell+2)}  $, that is,
$  \sigma=\cos( \bolds{\xi})=\cos ( \frac{\uppi}{2(\ell+2)}
) $.
In particular,
\begin{enumerate}
\item[]$  \bolds{\alpha}   =  0  $, if $  \sigma^2   =   1  /
2  $;
\item[]$  \bolds{\alpha}   =   -1   $, if $  \sigma^2   =   3
/  4  $; \ldots\ and
\item[]$  \bolds{\xi}   \downarrow  0  $, $  \bolds{\theta}_2
\downarrow  - \uppi  /   2  $, $  \bolds{\alpha}  \downarrow  -
\infty  $, as $  \sigma^2   \uparrow  1  $ in the limit.
\end{enumerate}

Note that the way of measuring the reflection angles in their paper is
to add
$  \uppi  /   2  $ to the angles $  (\bolds{\theta}_1, \bolds{\theta
}_2)   $
of Varadhan and Williams \cite{MR792398} but the parameter $  \bolds{\alpha}  $ is the same
in both
papers. Harrison \cite{MR509222}, Foschini \cite{MR651627} and Dai and Harrison \cite{MR1143393} also
studied the stationary distribution of the semimartingale reflected
Brownian motion.

 \noindent$\bullet $ In the case $  h > g >0  $, $  \sigma^2 =
\rho^2 = 1
/   2  $ with $  \bolds{\alpha} = 0 = \ell  $, $  \bolds{\theta}_1
= 0 =
\bolds{\theta}_2   $, we may compute the invariant distribution of the ranks
explicitly. From (\ref{eq: sum-of-exponential 1}), the stationary
density function $  \mathfrak p(\xi_1, \xi_2)   $ of $  (R_1(\cdot
), R_2(\cdot))  $ is given by\vspace*{-1pt}
\begin{equation}
\label{eq: sta dens 0} \mathfrak p(\xi_1, \xi_2) = 16
  h   (h - g)   \exp \bigl( - 4 (h   \xi_1 - g
\xi_2)   \bigr)  ,\qquad   0 < \xi_2 < \xi_1 <
\infty  .
\end{equation}
In fact, by direct computation the second term of (\ref{eq:
sum-of-exponential 1}) is zero. The first term of (\ref{eq:
sum-of-exponential 1}) is proportional to the exponential form $ \exp
( - 2
\langle\bolds{\mu}   , \mathbf{x} \rangle)   $. We obtain \eqref{eq:
sta dens 0} by
observing that $ \mathbf{x}=\sqrt{2 }  (\xi_1, \xi_2)' $ and $ \bm
\mu=\sqrt{2 }  (h,-g)' $. The value of the
normalizing constant comes from the fact that the invariant density of
$
(R_1(\cdot) - R_2(\cdot), R_2(\cdot))  $ is the product of
exponentials with
parameters $  (4 h ,  4 (h-g))  $.

\begin{rem}
The skew-symmetry condition of Harrison and Williams \cite{MR877593} holds for the process
$  (\sqrt{2 } R_1(\cdot), \sqrt{2 }  R_2(\cdot) )  $
in the case of equal variances. Under this skew-symmetry condition, the
invariant density has the form of a product of exponentials.

These parameters may be derived from the following heuristics. By (\ref
{2.15.a}) and (\ref{2.16}) and the strong law of large numbers for
Brownian motion, the local time grows linearly with
\begin{eqnarray*}
\lim_{t \to\infty} \bigl(L^{R_2}(t) / t \bigr) &=& \lim_{t\to\infty
}
\bigl(R_1(t) + R_2(t) - \xi- \nu  t - V(t) \bigr) / t =
h-g   ,
\\
 \lim_{t \to\infty} \bigl(L^{Y}(t) /t \bigr)&=& h   ,
\qquad   \mbox{a.s.}
\end{eqnarray*}
Thus, with these growth rates of the local times, the effective drift
rate of $  R_1(\cdot) - R_2(\cdot)  $ in (\ref{2.15.a}) is
heuristically $  - \lambda  t - \Lambda^{R_2}(t) \approx-(g + h)t
- (h - g)t = - 2  h  t   $ for the large $  t > 0  $. Similarly,
from (\ref{2.16}) the effective drift rate for $  R_2(\cdot)  $ is
$  (g - h) / (1/2) = - 2  (h - g)   $, where we divide $  (g - h)
$ by $  1  /  2  $ because the quadratic variation of $  \sigma
V_2(\cdot)  $ is a half of that of the standard Brownian motion.
Since the invariant density of Brownian motion in $  \mathbb R_+  $
with drift rate $  - \lambda< 0   $ reflected at the origin is known
to be exponential $  (2   \lambda)   $, we derive the parameters
$  (4   h, 4  ( h-g))  $ for $  (R_1(\cdot) - R_2(\cdot),
R_2(\cdot))  $; this is consistent with the consequence of (\ref{eq:
sta dens 0}) mentioned in the first paragraph of this Remark.
\end{rem}

\noindent$\bullet $ Similarly, in each case of $  \bolds{\alpha}   =   - \ell
  $, $
\ell\in\mathbb{N}  $ we may compute the invariant density of $
(R_1(\cdot),
R_2(\cdot))  $ from (\ref{eq: sum-of-exponential 1}). For example,
in the case
$  h > g > 0  $, $  \rho^2 = 1   /   4  $, $  \sigma^2 = 3
/   4
$ with
\[
  \bolds{\alpha}   =  -1 = - \ell  , \qquad    \bolds{\theta}_1 = 0  ,\qquad
    \bolds{\theta}_2 = - \uppi  /   6  ,  \qquad   \bolds{\xi} = \uppi
/   6  ,
\]
substituting $  \mathbf{x}'=(x_1, x_2) =  (2\xi_1,2\xi_2/\sqrt {3}   ) $ and
$  \bolds{\mu}'= (2h,-2g/\sqrt{3}   )$ into (\ref{eq:
sum-of-exponential 1}), we obtain for $  0 < \xi_2 < \xi_1 < \infty
  $ the invariant density $ \mathfrak p (\xi_1, \xi_2)  $ of $
(R_1(\cdot), R_2(\cdot))  $ as a linear combination of
\begin{eqnarray*}
\pi_0 (\mathbf{x})  &\propto&  \exp \bigl\{{- 8 [{h
\xi_1 - g\xi_2/3} ]} \bigr\} ,
\\
\exp \biggl\{{- \bolds{\mu}^\prime  %
\pmatrix{
\frac{3}{2} & - \frac{\sqrt{3}}{2}
\cr
\frac
{\sqrt{3}}{2} & \frac{3}{2} } %
\mathbf{x} }
\biggr\}   &=&   \exp \bigl\{{- \bigl[{(6h - 2 g) \xi_1 - 2 (g + h)
\xi_2 } \bigr]} \bigr\}    \quad   \mbox{and }
\\
\exp \biggl\{{- \bolds{\mu}^\prime  %
\pmatrix{
\frac{3}{2} & \frac{\sqrt{3}}{2}
\cr
\frac{\sqrt {3}}{2} & \frac{1}{2} } %
\mathbf{x}   }
\biggr\} & =& \exp \bigl\{{- \bigl[{(6h-2g)\xi_1+ (2h-2g/3)
\xi_2} \bigr]} \bigr\}   .
\end{eqnarray*}
Since $  \lim_{\sigma\uparrow1} \bolds{\alpha}  =   - \infty  $,
the invariant distribution of $ (R_1(\cdot), R_2(\cdot))  $ in
(\ref{2.20}), (\ref{2.21}) with $  \sigma=1 $ is conjectured to be
proportional to the infinite sum of exponentials as the limit of (\ref{eq:
sum-of-exponential 1}). It is an interesting open problem to determine
the invariant distribution for the degenerate case $  \sigma^2 = 1
$, as well as for general values of $  \sigma^2   $ that correspond
to noninteger scalars $  \bolds{\alpha}=( \bolds{\theta}_1+\bolds{\theta}_2) /
\bolds{\xi}=-\ell $.

 \begin{appendix}\label{app}

\section*{\texorpdfstring{Appendix: The other degenerate case, $\sigma=0$}{Appendix: The other degenerate case, sigma = 0}}
\label{A}
\renewcommand{\theequation}{A.\arabic{equation}}
\setcounter{equation}{0}


We have assumed throughout this work that the variance of the laggard
is positive. In this Appendix, we shall discuss briefly what happens
when the laggard undergoes a ``ballistic motion'' with positive drift
$  g >0 $, and the leader has unit variance, that is $  \rho=1 $
and $  \sigma=0 $.

Let us assume then that we have, on some filtered probability space $
( \Omega, \mathfrak{F}, \mathbb{P}),   \mathbf{F} = \{ \mathfrak
{F} (t) \}_{0
\le t < \infty} $, two continuous, nonnegative and adapted processes
$  X_1 (\cdot) $, $  X_2 (\cdot) $ that satisfy
\begin{eqnarray}
\label{A.1} \mathrm{d} X_1 (t) &= & ( g   \mathbf{
1}_{ \{ X_1 (t) \le
X_2 (t) \} } - h   \mathbf{ 1}_{ \{ X_1 (t) > X_2 (t) \} } ) \,  \mathrm{d} t  +
\mathbf{ 1}_{ \{ X_1 (t) > X_2 (t) \} } \, \mathrm {d} V (t),
\\
\label{A.2} \mathrm{d} X_2 (t) &=&  ( g   \mathbf{
1}_{ \{ X_1 (t) > X_2
(t) \} } - h   \mathbf{ 1}_{ \{ X_1 (t) \le X_2 (t) \} } )  \, \mathrm{d} t  +
\mathbf{ 1}_{ \{ X_1 (t) \le X_2 (t) \} } \, \mathrm{d} V (t) +  \mathrm{d} \Lambda(t)
\end{eqnarray}
for $  0 \le t < \infty $; here $  V(\cdot) $ is standard
Brownian motion, and $  \Lambda(\cdot) $ is a continuous, adapted
and nondecreasing process with $  \Lambda(0)=0 $ and
\begin{equation}
\label{A.3} \int_0^\infty\mathbf{
1}_{ \{ X_2 (t) >0 \} } \,  \mathrm{d} \Lambda (t)  = 0 , \qquad  \mbox{a.e.}
\end{equation}
The system of (\ref{A.1}), (\ref{A.2}) corresponds formally to that
of (\ref{3.11.b}), (\ref{3.12.b}), in light of the notation (\ref
{2.7}). No increasing component (such as
$  \Lambda(\cdot) $ of (\ref{A.2})) is needed in (\ref{A.1})
because, when the process $  X_1 (\cdot)  $ finds itself at the
origin, the positivity of its drift $  g >0 $ and the fact that its
motion is purely ballistic at that point are sufficient to ensure that
$  X_1 (\cdot)  $ stays nonnegative.

With the notation of (\ref{1.1}), (\ref{2.1}) it is fairly clear that
the difference $  Y (\cdot) = X_1 (\cdot) - X_2 (\cdot) $
satisfies the equation
\begin{equation}
\label{A.4} \mathrm{d} Y(t)  =  \sgn \bigl( Y(t) \bigr)   \bigl( - \lambda
\,\mathrm{d} t + \mathrm{d} V (t) \bigr) - \mathrm{d} \Lambda(t) ;
\end{equation}
we recall also the \textsc{Tanaka} formulae $  \mathrm{d} Y^+ (t) =
\mathbf{ 1}_{ \{ Y (t) > 0 \} } \,\mathrm{d} Y(t) + \mathrm{d} L^Y (t)
 $ and (\ref{Tanaka2}), the latter written now in the form
\begin{equation}
\label{A.6} \mathrm{d} \big| Y(t) \big|  =  - \lambda \, \mathrm{d} t + \mathrm{d} V
(t) - \sgn \bigl( Y(t) \bigr) \,  \mathrm{d} \Lambda(t) + 2  \, \mathrm{d}
L^Y(t) .
\end{equation}
With their help, we express the rankings of (\ref{ranks}) as
\begin{eqnarray}
\label{A.7} R_1 (t) &=&  X_2 (t) + Y^+ (t)  =
r_1 - h  t + V(t) + \int_0^t
\mathbf{ 1}_{ \{ X_1 (s) \le X_2 (s) \} }\,  \mathrm{d} \Lambda(s) + L^Y (t) ,
\\
\label{A.8} R_2 (t) &=&  X_1 (t) - Y^+ (t)  =
r_2 + g  t + \int_0^t \mathbf{
1}_{ \{ X_1 (s) > X_2 (s) \} }\,  \mathrm{d} \Lambda(s) - L^Y (t) .
\end{eqnarray}
We claim that we have the $  \mathbb{P}$-a.e. identities
\begin{equation}
\label{A.9} \int_0^\infty\mathbf{
1}_{ \{ X_1 (t) = X_2 (t) \} } \, \mathrm{d} t  = 0 ,\qquad   \int_0^\infty
\mathbf{ 1}_{ \{ X_1 (t) \neq X_2 (t)
\} }  \,\mathrm{d} \Lambda(t)  = 0 .
\end{equation}

Indeed, the first identity is a direct consequence of
\eqref{flat} and (\ref{A.6}). As for the second, we observe that for
any point $ t $ at which $  \mathbf{ 1}_{ \{ X_1 (\cdot) > X_2
(\cdot) \} } \, \mathrm{d} \Lambda(\cdot)  $ increases, we have $
X_2 (t)=0 $, thus $  R_2 (t)=0 $; but since $  g>0 $, we see from
(\ref{A.8}) that $ t $ must then be also a point of increase for $
L^Y (\cdot)  $, therefore $  X_1 (t) - X_2 (t) = Y(t)=0 $. We
conclude $  \int_0^\infty\mathbf{ 1}_{ \{ X_1 ( t) > X_2 ( t) \}
} \, \mathrm{d} \Lambda( t) \equiv0  $, thus
\[
  \Lambda( \cdot)  =  \int_0^{ \cdot} \mathbf{
1}_{ \{ X_1 (
t) \le X_2 ( t)=0 \} } \, \mathrm{d} \Lambda( t) =  \int_0^{
\cdot}
\mathbf{ 1}_{ \{ R_1 ( t) =0 \} }  \, \mathrm{d} \Lambda(t)
\]
in conjunction with (\ref{A.3}), and so
\[
\mathrm{d} ( R_1 - R_2 ) (t)  =  - \lambda\,
\mathrm{d} t + \mathrm{d} V (t) - \mathbf{ 1}_{ \{ R_1 ( t) =0 \} } \,  \mathrm {d}
\Lambda(t) + 2 \,  \mathrm{d} L^Y(t) ;
\]
a comparison with (\ref{A.6}) gives now the second identity of (\ref
{A.9}). In particular, we have shown that the process $  \Lambda
(\cdot) $ is supported on the set of visits by the process $  (X_1
(\cdot), X_2 (\cdot)) $ to the corner of the quadrant:
\begin{equation}
\label{A.10} \operatorname{supp} \bigl( \Lambda(\cdot) \bigr)   \subseteq  \bigl\{ t
\ge0\dvt R_1 (t) =0 \bigr\}  \subseteq  \bigl\{ t \ge0\dvt
R_2 (t) =0 \bigr\} .
\end{equation}

After all this, the equations of (\ref{A.7}), (\ref{A.8}) take the
particularly simple form
\begin{equation}
\label{A.11} R_1 (t)  =  r_1 - h  t + V(t) +
\Lambda(t) + L^Y (t) , \qquad  R_2 (t)  =
r_2 + g  t - L^Y (t) ,
\end{equation}
and give
\begin{eqnarray}
\label{A.5} X_1 (t) + X_2 (t)  &=&  R_1
(t) + R_2 (t)  =  \xi+ \nu  t + V(t) + \Lambda(t)  ,
\\
\label{A.12} R_1 (t) - R_2 (t)  &=&  |y| - \lambda  t
+ V(t) + \Lambda(t) + 2  L^Y (t) .
\end{eqnarray}

\noindent$\bullet $ We observe from (\ref{A.10}), (\ref{flat}) that $
\operatorname{supp}  ( \Lambda(\cdot) + 2 L^Y (\cdot)  )
\subseteq \{ t \ge0\dvt R_1 (t) - R_2 (t)=0 \}  $, so it
follows from (\ref{A.12}) that the process $  R_1 (\cdot) - R_2
(\cdot) \ge0 $ is the \textsc{Skorokhod} reflection at the origin
of the Brownian motion with negative drift
\begin{equation}
\label{A.16} Z(t) =  |y| - \lambda  t + V(t) ,\qquad   0 \le t < \infty ,
\end{equation}
namely
\begin{equation}
\label{A.14} \Lambda(t) + 2  L^Y (t)  =  \max_{0 \le s \le t }
\bigl( - Z(s) \bigr)^+  =  \max_{0 \le s \le t }  \bigl( - |y| + \lambda  s -
V(s) \bigr)^+ ,\qquad   0 \le t < \infty .
\end{equation}

\noindent$\bullet $ On the other hand, we observe from (\ref{A.11}) that
\[
0 \le2  R_2 (t) =  2  r_2 + 2 g  t - \bigl(
\Lambda(t) + 2  L^Y (t) \bigr) + \Lambda(t) ,\qquad   0 \le t < \infty
 ,
\]
so in light of (\ref{A.10}) and the theory of the \textsc{Skorokhod}
reflection problem once again, we obtain
\begin{equation}
\label{A.15} \Lambda(t)  =  \max_{0 \le s \le t }  \bigl( - 2
r_2 - 2  g s + \Lambda(s) + 2  L^Y (s) \bigr)^+ ,
\qquad  0 \le t < \infty .
\end{equation}

\renewcommand{\therem}{A.\arabic{rem}}
\setcounter{rem}{0}

\begin{rem}
Likewise, from (\ref{A.10}) the support of $  \Lambda(\cdot) $ is
included in the zero-set of the process $  R_1 (\cdot) + R_2 (\cdot)
\ge0 $; then (\ref{A.5}) and the theory of the \textsc{Skorokhod}
reflection problem give
\begin{equation}
\label{A.13} \Lambda(t)  =  \max_{0 \le s \le t }  \bigl( - \xi- \nu s - V
(s) \bigr)^+ ,\qquad   0 \le t < \infty .
\end{equation}
In other words, the sum $  X_1 (\cdot) + X_2 (\cdot)= R_1 (\cdot) +
R_2 (\cdot)  $ is Brownian motion with drift $  \nu= g-h $ and
reflection at the origin;
we are indebted to \textsc{Dr. Phillip Whitman} for this observation.

Consequently, if $  h \ge g $, the process $  (X_1 (\cdot), X_2
(\cdot)) $ visits the corner of the nonnegative quadrant with
probability one; whereas, if $  h <g $, we have
\[
\mathbb{P}  \bigl(  X_1 (t) = X_2 (t) =0 ,
\mbox{for some }    t \ge0  \bigr)  =  \mathrm{e}^{  - 2  (g-h)  \xi} .
\]
\end{rem}

\begin{rem}
Applying the \textsc{Tanaka} formula (\ref{Tanaka2}) to the
continuous, nonnegative semimartingales $  X_1 (\cdot) + X_2 (\cdot)
 $ in (\ref{A.5}) and $  R_1 (\cdot) - R_2 (\cdot) $ in (\ref
{A.12}), we obtain the identifications
\begin{equation}
\label{A.17} \Lambda(\cdot) =  L^{X_1 +X_2} (\cdot) ,\qquad   \Lambda(
\cdot) + 2  L^Y (\cdot)  =  L^{R_1 -R_2} (\cdot) ,
\end{equation}
thus also
\[
L^{|X_1 -X_2|} (\cdot) - 2  L^{X_1 -X_2} (\cdot) =  L^{X_1 +X_2} (
\cdot) .
\]
On the other hand, we have the $  \mathbb{P}$-a.e. properties $
\int_0^\infty\mathbf{ 1}_{\{R_1(t)=0\}}  \, \mathrm{d} t =0$, $ \int_0^\infty\mathbf{ 1}_{\{R_1(t)=R_2 (t)\}}  \, \mathrm{d} t =0  $
from (\ref{A.11}), (\ref{A.12}) and (\ref{flat}). In conjunction
with (\ref{2.19.a}) and (\ref{Tanaka}) -- in particular, the fact
that $  L^X(\cdot) \equiv0 $ holds for a continuous semimartingale
$ X(\cdot) $ of finite variation -- we obtain from these equations
and (\ref{A.10}) the identifications
\begin{eqnarray*}
L^{R_1} (\cdot)  &=&  \int_0^\cdot
\mathbf{ 1}_{ \{ R_1 (t)=0 \} } \, \mathrm{d} R_1(t)  =  \Lambda(\cdot)
+ \int_0^\cdot\mathbf{ 1}_{ \{ R_1 (t)=0 \} }
\,\mathrm{d} L^Y (t) ,
\\
0 =  L^{R_2} (\cdot)  &= & \int_0^\cdot
\mathbf{ 1}_{ \{ R_2
(t)=0 \} } \, \mathrm{d} R_2(t)  =  g \int
_0^\cdot\mathbf{ 1}_{ \{
R_2 (t)=0 \} } \, \mathrm{d} t
- \int_0^\cdot\mathbf{ 1}_{ \{ R_2
(t)=0 \} }
\,\mathrm{d} L^Y (t) ,
\end{eqnarray*}
thus also
\begin{equation}
\label{A.18} \int_0^\cdot\mathbf{
1}_{ \{ R_1 (t)=0 \} } \, \mathrm{d} L^Y (t)  =  g \int
_0^\cdot\mathbf{ 1}_{ \{ R_1 (t)=0 \} } \, \mathrm{d} t
  =  0 ,\qquad   \Lambda(\cdot) =  L^{X_1 \vee X_2} (\cdot) .
\end{equation}
It is rather interesting that the same process $  \Lambda(\cdot) $ should
do ``triple duty'', as the local time of both the sum $  X_1 (\cdot)
+ X_2
(\cdot) $ and the maximum $  X_1 (\cdot) \vee X_2 (\cdot) $, and
as the
increasing process in the \textsc{Skorokhod} reflection for the
minimum $
X_1 (\cdot) \wedge X_2 (\cdot) $.
\end{rem}

\textit{Synthesis}: Now we can reverse the above steps. Starting with a
standard Brownian motion $  V(\cdot) $, we define $  Z (\cdot) $,
$  \Lambda(\cdot) $ and $  L^Y (\cdot) $ via (\ref
{A.16})--(\ref{A.15}); then $  R_1 (\cdot) $, $ R_2 (\cdot) $ via
(\ref{A.11}), and
\begin{equation}
\label{A.19} G(\cdot)  :=  R_1 (\cdot)- R_2 (\cdot)
 =  Z(\cdot) + \Lambda (\cdot) + 2  L^Y(\cdot)
\end{equation}
as in (\ref{A.12}). It is clear from (\ref{A.14}) that this process
$  G(\cdot) $ is the \textsc{Skorokhod} reflection at the origin of
the Brownian motion with negative drift $  Z(\cdot) $ in (\ref
{A.15}), thus nonnegative. It is also clear that the process $
\Lambda(\cdot) $ satisfies (\ref{A.10}), as well as (\ref
{A.13})--(\ref{A.18}).

On a suitable extension of the probability space we need to find now a
continuous semimartingale $  Y(\cdot) $ of the form (\ref{A.4}),
with the help of which we can ``unfold'' the process $  G(\cdot) $
of (\ref{A.19}) in the form $  G(\cdot) = | Y(\cdot)| $. Once this
has been done we can define
\[
X_1 (\cdot)  :=  R_2 (\cdot) + Y^+ (\cdot) ,\qquad
X_2 (\cdot)  :=  R_1 (\cdot) - Y^+ (\cdot)
\]
and verify the equations (\ref{A.1}), (\ref{A.2}) in a
straightforward manner. In order to carry out this unfolding, the
method outlined in Section \ref{sec4.1.a} is inadequate; it has to
be modified as follows.

We enumerate the excursions of $  G(\cdot) $ away from the origin,
just as before, but now distinguish between those that originate at the
corner of the quadrant ($R_1 =0$), and the rest ($R_1 >0$). Excursions
of the first type are always marked $ \Phi=-1 $; while excursions of
the second type are assigned marks $  \Phi= \pm1 $ independently of
each other, and with equal probabilities ($1/2$), just as in Section
\ref{sec4.1.a}. The resulting process $  Y(\cdot)= \Phi(\cdot)
G(\cdot)  $ satisfies $  G(\cdot) = | Y(\cdot)| $ and
\begin{eqnarray*}
\mathrm{d} Y(t)  &= & \overline{\sgn} \bigl( Y(t) \bigr)   \mathbf{
1}_{ \{ R_1 (t)
>0 \} } \, \mathrm{d} G(t) - \mathbf{ 1}_{ \{ R_1 (t)=0 \} } \, \mathrm {d}
G(t)
\\
 &=&  \overline{\sgn} \bigl( Y(t) \bigr)   \mathbf{ 1}_{ \{ R_1
(t) >0 \} } \cdot
\mathbf{ 1}_{ \{ G (t) >0 \} }  \, \mathrm{d} G(t) - \mathbf{ 1}_{ \{ R_1 (t)=0 \} } \cdot
\mathbf{ 1}_{ \{ R_1 (t)=R_2
(t) \} } \, \mathrm{d} ( R_1 - R_2 )
(t)
\\
 &=&  \overline{\sgn} \bigl( Y(t) \bigr)   \mathbf{ 1}_{ \{ R_1
(t) >0 \} } \cdot
\mathbf{ 1}_{ \{ G (t) >0 \} } \,  \mathrm{d} Z(t) - \mathbf{ 1}_{ \{ R_1 (t)=0 \} }
\,\mathrm{d} L^{ R_1 - R_2} (t)
\\
 &= & \overline{\sgn} \bigl( Y(t) \bigr)   \mathbf{ 1}_{ \{ G (t)
>0 \} }
\,\mathrm{d} Z(t) - \mathbf{ 1}_{ \{ R_1 (t)=0 \} }  \, \mathrm{d} L^{ R_1 - R_2} (t)
\\
 &= & \overline{\sgn} \bigl( Y(t) \bigr)   \bigl( - \lambda \,\mathrm{d} t +
\mathrm{d} V(t) \bigr) - \mathbf{ 1}_{ \{ R_1 (t)=0
\} }   \bigl( \mathrm{d}
\Lambda(t) + 2 \, \mathrm{d} L^Y (t) \bigr)
\\
 &=&  \sgn \bigl( Y(t) \bigr)   \bigl( - \lambda\,\mathrm{d} t + \mathrm{d} V(t)
\bigr) - \mathrm{d} \Lambda(t) ,
\end{eqnarray*}
that is, the equation of (\ref{A.4}), as promised. We have used (\ref
{2.19.a}), (\ref{A.10}), (\ref{A.17}), (\ref{A.18}) as well as the
$  \mathbb{P}$-a.e. properties $ \int_0^\infty\mathbf{ 1}_{\{
G(t)=0\}}
\,\mathrm{d} t =0 $, $  \int_0^\infty\mathbf{ 1}_{\{G(t)>0\}}
\,\mathrm{d} (G-Z)(t) =0 $; the first of these is a consequence of
(\ref{flat}), and the second of \textsc{Skorokhod} reflection.
\end{appendix}

\section*{Acknowledgements}
We are grateful to Dr. E.$\,$Robert Fernholz for posing this problem
    in his report Fernholz \cite{Fernholz2011a},  for encouraging us to pursue it, for giving us permission to reproduce here Figure 1 from his report, and for
    supplying an elegant and highly original argument for the proof of
    Proposition~\ref{Prop3}. We also thank Professor Ruth Williams for her
    interest and her very helpful suggestions,   Dr. Johannes Ruf for his many
    comments on successive versions of this work,   Dr. Phillip Whitman for
    discussions that helped us sharpen our understanding, and the referees for their careful readings and helpful suggestions.
Ioannis Karatzas research was supported in part by National Science Foundation Grant DMS-09-05754.
Vilmos Prokaj research was supported    by the European Union and co-financed by the European Social Fund (grant agreement no. TAMOP 4.2.1/B-09/1/KMR-2010-003).


%

\printhistory

\end{document}